\documentclass[11pt,reqno]{amsart}
\usepackage[a4paper,margin=28mm]{geometry}
\usepackage{amsmath,amssymb,mathtools,booktabs,longtable,array}
\usepackage{enumitem,microtype,needspace}
\usepackage[T1]{fontenc}
\usepackage{lmodern}
\PassOptionsToPackage{hyphens}{url}
\usepackage[colorlinks=true,linkcolor=blue,citecolor=blue,urlcolor=blue]{hyperref}
\numberwithin{equation}{section}
\newtheorem{theorem}{Theorem}[section]
\newtheorem{lemma}[theorem]{Lemma}
\newtheorem{proposition}[theorem]{Proposition}
\newtheorem{corollary}[theorem]{Corollary}
\newtheorem*{conjecture}{Conjecture}
\theoremstyle{definition}

\newtheorem{remark}[theorem]{Remark}

\newcommand{\Z}{\mathbb Z}
\newcommand{\Q}{\mathbb Q}
\newcommand{\R}{\mathbb R}
\newcommand{\C}{\mathbb C}
\newcommand{\HH}{\mathbb H}
\newcommand{\II}{\mathrm{II}}
\newcommand{\ee}{\mathfrak e}
\newcommand{\rk}{\operatorname{rk}}
\newcommand{\Num}{\operatorname{Num}}
\newcommand{\Pic}{\operatorname{Pic}}
\newcommand{\Span}{\operatorname{span}}
\DeclareMathOperator{\SL}{SL}
\setlist[enumerate]{label=\textup{(\roman*)},leftmargin=*}
\allowdisplaybreaks[2]

\title[The Borisov--Nuer conjecture and root decompositions]
{The Borisov--Nuer conjecture and root decompositions in even unimodular Lorentzian lattices}
\author{Haowu Wang}
\address{School of Mathematics and Statistics, Wuhan University, Wuhan 430072, Hubei, China}
\email{haowu.wangmath@whu.edu.cn}
\author{Shaoyun Yi}
\address{School of Mathematical Sciences, Xiamen University, Xiamen, Fujian 361005, China}
\email{yishaoyun926@xmu.edu.cn}
\date{\today}

\begin{document}
\begin{abstract}
We prove that every vector of $\mathrm{II}_{8n+1,1}$, with $n\ge1$,
is a sum of two roots if and only if $n=1$ or $2$.
This settles the Borisov--Nuer conjecture and implies the existence
of Ulrich line bundles for every very ample polarization on a complex
unnodal Enriques surface. For vectors of negative norm, we express
the number of root decompositions as a Fourier coefficient of a
vector-valued theta series, reducing existence to positivity.
The main analytic difficulty is to obtain estimates uniform over
lattices of unbounded determinant. We address this by combining
a reduction based on primitive sublattices of codimension one and
least determinant with explicit bounds for theta norms and Fourier
coefficients of cusp forms.
\end{abstract}

\begingroup
\raggedbottom
\maketitle
\tableofcontents
\clearpage
\endgroup

\section{Introduction and main results}\label{sec:introduction}

Let $X$ be a smooth complex projective variety of dimension $d$, and let
$H$ be a very ample divisor. A vector bundle $\mathcal E$ on $X$ is
\emph{Ulrich with respect to $H$} if
\[
 H^i(X,\mathcal E(-jH))=0
 \qquad\text{for all }i\text{ and }1\le j\le d.
\]
Ulrich bundles provide linear resolutions of the corresponding sheaves
in projective space and enter the study of Chow forms and cohomology
tables; see \cite{ES03,BeauvilleSurvey}. Their existence is a basic
problem, and the rank-one case imposes particularly strong restrictions
on the divisor classes of the variety.

Let $Y$ be a complex Enriques surface. We equip its numerical divisor
lattice with the negative of the intersection form. For lattices,
$\oplus$ denotes the orthogonal direct sum; for abelian groups, it denotes
the ordinary direct sum. Thus
\[
 \Num(Y)=\Pic(Y)/\langle K_Y\rangle\cong U\oplus E_8\cong\II_{9,1},
 \qquad (D,E)=-D\mathbin{\cdot}E,
\]
where $\cdot$ denotes the usual intersection product; see \cite{CD89}
and \cite[\S3]{Dolgachev2016}.
Here $U$ is a hyperbolic plane and $E_8$ is the positive definite
root lattice. We write $D^2=(D,D)$, so a very ample class $H$ has $H^2<0$.
This convention, with one negative direction, will be used throughout.
The surface $Y$ is \emph{unnodal} if it contains no smooth rational
curves. Borisov and Nuer \cite[Proposition~2.1]{BN18} proved that
a line bundle $\mathcal O_Y(D)$ on such a surface has vanishing
cohomology in every degree if and only if $D^2=2$.
Consequently, $\mathcal O_Y(D)$ is $H$-Ulrich if and only if
\[
                  (D-H)^2=(D-2H)^2=2.
\]
Setting $r=D-H$ and $s=2H-D$, we obtain the equivalent condition
\begin{equation}\label{intro:BN}
                         H=r+s,\qquad r^2=s^2=2.
\end{equation}
Borisov and Nuer established this decomposition for multiples of a
Fano polarization and proposed the following numerical conjecture
\cite[Conjecture 2.2 and Theorem 2.4]{BN18}, stated here in our sign
convention.

\begin{conjecture}[Borisov--Nuer]
Every vector of $\II_{9,1}$ is a sum of two vectors of square $2$. 
\end{conjecture}

Aprodu and Kim \cite{AK20} related this numerical condition to the
images of moduli spaces of polarized Enriques surfaces in K3 moduli
spaces. Their Conjecture~1.3 is a direct consequence of Borisov--Nuer's conjecture.  

We prove the conjecture and study the same question for all even
unimodular Lorentzian lattices. For an integer $n\ge1$, let
$\II_{8n+1,1}$ denote the
even unimodular lattice of signature $(8n+1,1)$, and put
\[
 \mathcal R_n=\{r\in\II_{8n+1,1}:r^2=2\}
\]
for the set of roots. A lattice is called rootless if it has no roots.
We define
\begin{equation}\label{intro:count}
 R_n(\gamma)=\#\{(r,s)\in\mathcal R_n^2:r+s=\gamma\},
 \qquad\gamma\in\II_{8n+1,1}.
\end{equation}
Representations are ordered, and we write $R_n(\gamma)=\infty$ when
there are infinitely many.

The integral orthogonal group gives another formulation of the
conjecture. Choose a hyperbolic plane with $e^2=f^2=0$ and $(e,f)=-1$.
Since $\mathrm{O}(\II_{9,1})$ acts transitively on roots
by Lemma~\ref{gen:lem:positive-normal}, the conjecture is equivalent to
\begin{equation}\label{intro:group-problem}
 \II_{9,1}=\mathrm{O}(\II_{9,1})(e-f)+\mathrm{O}(\II_{9,1})(e-f).
\end{equation}
The following theorem determines which of these lattices have the
decomposition property, describes the vectors of nonnegative norm,
and classifies the obstructions of sufficiently large negative norm.

\begin{theorem}\label{thm:classification}
Let $n\ge1$ and $\gamma\in\II_{8n+1,1}$.
\begin{enumerate}
\item Every vector of $\II_{8n+1,1}$ is a sum of two roots if and only if
$n=1$ or $n=2$.
\item If $\gamma^2<0$, then $R_n(\gamma)$ is finite.
If $\gamma^2>0$ or $\gamma=0$, then $R_n(\gamma)=\infty$.
\item If $\gamma^2=0$ and $\gamma\ne0$, write $\gamma=ky$, where $k\ge1$ and $y\in \II_{8n+1,1}$ is primitive isotropic. Then $P_y=(y^\perp\cap \II_{8n+1,1})/\Z y$ is an even positive-definite unimodular lattice of rank $8n$. Moreover, 
\[
R_n(\gamma)=
\begin{cases}
0,&P_y\text{ has no roots},\\
\infty,&P_y\text{ has a root}.
\end{cases}
\]
\item If $n\ge3$, then for any even integer $T>0$, there is
a primitive $\gamma\in \II_{8n+1,1}$ with $\gamma^2=-T$ and
$R_n(\gamma)=0$. There are also primitive isotropic vectors
without representation.
\item There is $T_n^*>0$, depending only on $n$, such that, if
$-\gamma^2\ge T_n^*$, then $R_n(\gamma)=0$ if and only if
$\gamma$ is primitive and $\gamma^\perp\cap\II_{8n+1,1}$ contains
a rootless even unimodular sublattice of rank $8n$.
\end{enumerate}
\end{theorem}

Thus the isotropic case in higher rank depends on the cusp lattice.
The proof also gives sharp quantitative statements in the two dimensions
where every vector is represented. We write $\gamma\sim\gamma'$ when
the vectors lie in the same orbit under the full integral orthogonal group.

\begin{corollary}\label{cor:ii9-sharp}
Let $\gamma\in\II_{9,1}$ with $\gamma^2=-T<0$.
\begin{enumerate}
\item $R_1(\gamma)\ge480$. Equality holds if and only if $\gamma$ is
primitive and $\gamma^\perp\cap \II_{9,1}$ contains $E_8$. Equivalently,
in some decomposition $\II_{9,1}=U\oplus E_8$,
\[
 \gamma=a_0e+b_0f,\qquad a_0,b_0>0,\quad\gcd(a_0,b_0)=1.
\]
\item If $R_1(\gamma)\ne480$, then $R_1(\gamma)\ge672$.
Equality holds if and only if $T=4$ and
$\gamma^\perp\cap \II_{9,1}\cong D_9$; these vectors form one orbit.
\item If $\gamma$ is nonprimitive, then $R_1(\gamma)\ge2640$.
Equality holds if and only if $\gamma\sim2(e+f)$.
\end{enumerate}
\end{corollary}

For the next statement, let $\mathcal U$ be the set of primitive
negative-norm vectors $\gamma$ in $\II_{17,1}$ whose perpendicular contains an even
positive-definite unimodular lattice of rank $16$. Equivalently, in some decomposition $\II_{17,1}=U\oplus P$, 
\[
\gamma=a_0e+b_0f,\qquad a_0,b_0>0,\quad\gcd(a_0,b_0)=1,
\]
where $P\cong E_8\oplus E_8$ or $D_{16}^{+}$. Here $D_{16}^{+}$ denotes the even unimodular extension of index two of the root lattice $D_{16}$. 

\begin{corollary}\label{cor:ii17-sharp}
Let $\gamma\in\II_{17,1}$ with $\gamma^2=-T<0$.
\begin{enumerate}
\item $R_2(\gamma)\ge960$, with equality precisely for
$\gamma\in\mathcal U$.
\item If $\gamma\notin\mathcal U$, then $R_2(\gamma)\ge1632$.
Equality holds in exactly two orbits, both with $T=2$; their
perpendiculars are the even overlattices of $A_{17}$ and
$E_7\oplus D_{10}$ of indices three and two, respectively.
\item If $\gamma\notin\mathcal U$ and $T\ge4$, then
$R_2(\gamma)\ge4488$. Equality holds in one orbit with $T=4$,
whose perpendicular is the even overlattice of $E_6\oplus A_{11}$
of index three.
\item If $\gamma$ is nonprimitive, then $R_2(\gamma)\ge55488$.
Equality holds precisely for the doubles of the two vectors in
\textup{(ii)}, up to isometry.
\end{enumerate}
\end{corollary}

\subsection{Outline of the proof}
The main work concerns vectors of negative norm. We obtain coefficient
estimates that are uniform as the perpendicular lattice varies.
Lattice reduction then leaves finitely many vectors to examine for $n=1,2$.
The estimates also apply in higher rank; there, rootless unimodular
summands account for all exceptions of sufficiently large negative norm.
We explain the argument in six steps. Put $T=-\gamma^2$ and write
each nonzero vector as $\gamma=ky$, with $k\ge1$ and $y$ primitive.

\vspace{2mm}

\emph{Step 1: vectors of nonnegative norm.}
For $T<0$, the primitive vector $y$ is isometric to $e-df$ with $d>0$.
The fact that $E_8$ represents every positive even integer supplies
infinitely many root pairs, by varying their $e$-coordinates. 
For $T=0$ and $\gamma\ne0$, identify $P_y$ with the orthogonal
complement of a split hyperbolic plane containing $y$. Then
\[
 (r,s)=(ay+z,(k-a)y-z),\qquad a\in\Z,\quad z\in P_y,
 \quad z^2=2.
\]
Thus, the answer is $0$ or $\infty$ depending on whether the cusp lattice
$P_y$ is rootless or contains a root. In ranks $8$ and $16$, it
always contains roots. The zero vector also has infinitely many root
decompositions. 

\vspace{2mm}

\emph{Step 2: representations as theta coefficients.}
For $T>0$, put $y^2=-\delta$ and
$K=y^\perp\cap\II_{8n+1,1}$, positive definite of rank $8n+1$.
Unimodular gluing gives
$\det K=\delta$, $A_K=K'/K\cong\Z/\delta\Z$, and a class
$\mu\in A_K$ of order at most two. Since $T=k^2\delta$, we have
$\delta\le T$. Lemma~\ref{r9:lem:gluing-shell} gives
\begin{equation}\label{intro:shell}
 R_n(\gamma)
 =\#\{v\in K+\mu:Q(v)=T/8+1\}
 =a_{K}(\mu,T/8+1).
\end{equation}
Here $Q(v)=v^2/2$, $K+\mu$ is a coset in $K'$, and
$a_{K}(\mu,t)$ counts the vectors $v\in K+\mu$ with $Q(v)=t$.
The bijection is $r=\gamma/2+v$, $s=\gamma/2-v$.
In particular, the count is finite, and existence is equivalent to
positivity of this theta coefficient. As $\gamma$ varies, so do $K$,
its determinant and its Weil representation. An asymptotic estimate
with a constant depending on $K$ would therefore not give the
uniform bound on $T$ needed for a finite reduction.

The main step is to compare $K$ with a suitable sublattice.
Let $c$ be the least determinant of a primitive sublattice
$S\subset K$ of rank $8n$. Then
\[
 \lambda_1(K')^2=c/\delta,
\]
where $\lambda_1(K')$ is the length of a shortest nonzero vector of $K'$.
Choose $S=u^\perp\cap K$ with $u\in K'$ attaining this minimum. Projecting $K+\mu$ onto $S_\R^\perp$ and taking the closest fiber gives, by
Lemma~\ref{r9:lem:fibre},
\begin{equation}\label{intro:fibre}
 R_n(\gamma)\ge a_{S}(\mu_S,t),\qquad
 t\ge\frac T8\left(1-\frac1c\right)+1.
\end{equation}
We estimate the coefficient of $\Theta_K$ when $c$ is large and that of $\Theta_S$
when $c$ is small. More precisely, for a fixed $b>0$,
\[
 c\ge b\sqrt T\ \Longrightarrow\
 \lambda_1(K')^2\ge b/\sqrt T,\qquad
 c\le b\sqrt T\ \Longrightarrow\ \det S\le b\sqrt T.
\]
The first condition controls the short vectors of $K'$; the second
gives a determinant bound for the fiber lattice. The fiber inequality
holds for every primitive corank-one $S$. The minimizing choice is
needed for the first implication, through
$\lambda_1(K')^2=c/\delta$.
\vspace{2mm}

\emph{Step 3: unimodular sublattices.}
When $c=1$, a scalar theta series gives an exact formula.
Indeed, $\II_{8n+1,1}=U\oplus P$ with $P$ unimodular and
$\gamma=k(a_0e+b_0f)$, $a_0,b_0>0$, $\gcd(a_0,b_0)=1$.
For $n=1,2$, the scalar theta identities in ranks $8$ and $16$ give
\begin{equation}\label{intro:split9}
 R_n(\gamma)
   =240n\sum_{j=0}^k\sigma_{4n-1}\bigl(1+a_0b_0j(k-j)\bigr),
\end{equation}
where $\sigma_j(m)=\sum_{d\mid m}d^j$; see
Corollary~\ref{gen:cor:small-split-counts}. The minimum counts in the
case $c=1$ are therefore $480$ and $960$, respectively.

\vspace{2mm}

\emph{Step 4: comparison of Eisenstein and cusp coefficients.}
Suppose $c>1$. Then $c\ge4$, so the fiber index in
\eqref{intro:fibre} is at least $3T/32+1$.
For $L=K$ or $S$, put $r=\rk L$, $D=\det L$ and write
$\Theta_L=E_L+F_L$, where $E_L$ is the
Eisenstein series with constant term $\ee_0$ and $F_L$ is a cusp form.
Since $A_K$ is cyclic and $A_S$ has at most two generators,
the local densities have uniform positive lower bounds, giving
\begin{equation}\label{intro:main-term}
 a_{E_L}(\mu,t)\ge
 \eta_r\frac{(2\pi)^{r/2}}{\Gamma(r/2)}
       \frac{t^{r/2-1}}{\sqrt D}
\end{equation}
at the relevant compatible indices, with explicit $\eta_r>0$
depending only on $r$.

The Poincar\'e-series bound in Lemma~\ref{r9:lem:poincare-bound}
controls the cusp term through its Petersson norm, with
discriminant-group factors bounded using the same restrictions on
$A_K$ and $A_S$. Complementary determinants bound products of
successive minima of $L'$; lattice-point counting and integration
then give uniform theta norm estimates. For $n=1,2$, these yield
\begin{equation}\label{intro:norm-bounds}
\begin{array}{c|cc}
 &\|\Theta_K-\ee_0\|^2< &\|\Theta_S-\ee_0\|^2<\\\hline
 n=1&0.017\,T^{7/4}&0.014\,T^{3/2}\\
 n=2&0.0015\,T^{15/4}&0.0013\,T^{7/2}.
\end{array}
\end{equation}
The columns apply when $c\ge b\sqrt T$ and $c\le b\sqrt T$,
respectively, for $T\ge T_0$, with
$(b,T_0)=(8/5,200)$ for $n=1$ and $(3/2,100)$ for $n=2$;
see Corollary~\ref{norm:cor:thick-thin}.

Combining these estimates with
$\|F_L\|\le\|\Theta_L-\ee_0\|+1/\sqrt{24}$ gives
\[|a_{F_L}(\mu,t)|<(19/20)a_{E_L}(\mu,t)\]
at the selected indices in each branch, for $T\ge T_0$.
Proposition~\ref{gen:prop:analytic-range} therefore yields
\begin{equation}\label{intro:large9}
 R_1(\gamma)>\frac{T^{11/4}}{200}>2640\quad(T\ge200),\qquad
 R_2(\gamma)>\frac{T^{27/4}}{10^6}>62880\quad(T\ge100).
\end{equation}
Hence equality for $c>1$ is confined to bounded $T$.

\vspace{2mm}

\emph{Step 5: bounded ranges and equality cases.}
The reduction bounds in Lemmas~\ref{r9:lem:chamber},
\ref{r9:lem:enumeration}, \ref{r17:lem:reduction}
and~\ref{r17:lem:finitebound} give finite lists meeting every
isometry orbit with $0<T\le200$ in $\II_{9,1}$ and
$0<T\le100$ in $\II_{17,1}$.

For each reduced $\gamma=ae+bf+z$ with $b>0$, write
$r=pe+qf+v$, where $q\in\Z$ and $z,v\in E_8^{\oplus n}$,
$n=1,2$. The root equations become
\[
 \lVert bv-qz\rVert^2=2b^2+Tq(b-q),\qquad
 p=\frac{T/2-aq+(z,v)}b\in\Z.
\]
Nonnegativity bounds $q$ to a finite interval, and each $q$
leaves finitely many lattice points on a sphere.
The resulting exact counts and lower bounds are recorded in
Lemmas~\ref{comp:ii9} and~\ref{comp:ii17}.

Together with \emph{Steps 3} and \emph{4}, these tables leave only
the split equality cases and finitely many exceptional candidates.
Lemmas~\ref{r9:lem:D9-count} and~\ref{r17:lem:exceptional-shells}
identify the exceptional perpendicular lattices and verify the
exact counts, completing Corollaries~\ref{cor:ii9-sharp}
and~\ref{cor:ii17-sharp}.

\vspace{2mm}

\emph{Step 6: higher ranks.}
For every $n\ge3$, choose a rootless even positive-definite
unimodular lattice $P_n$ of rank $8n$ and identify
$\II_{8n+1,1}=U\oplus P_n$.
The cusp at $e$ has no roots. For a root decomposition $r+s$ of $e+df$,
$d\ge1$, the root equations force $r=pe+d(1-p)f+z$ and
\[
 z^2=2+2dp(1-p),\qquad  p\in\Z,\quad z\in P_n. 
\]
The right side is negative, unless $p=0,1$, when it requires a root.
This gives the isotropic and negative-norm counterexamples.
The coefficient argument for $c>1$ still applies. For each fixed
$n$, it gives
\[
 R_n(\gamma)\gg_n T^{4n-5/4}
\]
for sufficiently large $T$, with a bound depending only on $n$;
see Proposition~\ref{gen:prop:analytic-range}.
When $c=1$, the scalar theta asymptotic excludes sufficiently large
nonprimitive exceptions. The remaining exceptions are precisely the primitive
vectors whose perpendicular contains a rootless unimodular summand
of rank $8n$, as asserted in Theorem~\ref{thm:classification} \textup{(v)}.

\vspace{2mm}

The remainder of this paper is organized as follows. Section~\ref{sec:lattice-reductions} gives the common lattice reductions. Section~\ref{sec:theta-estimates} establishes the uniform estimates for theta coefficient. Sections~\ref{sec:ii9} and \ref{sec:ii17} prove the sharp counts for $n=1$ and $n=2$, respectively. Section~\ref{sec:higher-rank} treats the higher-rank obstructions and the diameters of the Cayley graph, and completes the proof of Theorem~\ref{thm:classification}. In particular, we show that every vector in $\II_{8n+1,1}$ is a sum of at most three roots. The finite enumerations are included in Sections~\ref{sec:ii9} and~\ref{sec:ii17}.

\section{Lattice reductions}\label{sec:lattice-reductions}
\subsection{Notations and preliminaries}

An integral lattice $L$ is a free abelian group of finite rank with a
nondegenerate integral symmetric bilinear form $(\ ,\ )$.
Write $x^2=(x,x)$ and $Q(x)=x^2/2$; $L$ is even if $Q(L)\subset\Z$.
Its determinant is the absolute determinant of a Gram matrix.
For an even lattice, set
\[
 L'=\{x\in L\otimes\Q:(x,L)\subset\Z\},\qquad
 A_L=L'/L,\qquad q_L(x+L)=Q(x)+\Z.
\]
Thus $A_L$ is the discriminant group, $q_L$ its discriminant
quadratic form, and $|A_L|=\det L$; $L$ is unimodular if
$\det L=1$. For $x\in L'$, write $[x]=x+L\in A_L$.
We write $L_\Q=L\otimes_\Z\Q$ and $L_\R=L\otimes_\Z\R$.
For a finite abelian group $A$, let $\ell(A)$ be its
least number of generators. A sublattice $M\subset L$ is primitive
if $L/M$ is torsion-free; a nonzero $y\in L$ is primitive if
$L\cap\Q y=\Z y$.

The first entry of a signature is positive, and $n\ge1$ throughout
the lattice reductions. Fix
\[
 U=\Z e+\Z f,\qquad e^2=f^2=0,\qquad(e,f)=-1.
\]
The classification of indefinite even unimodular lattices gives
$\II_{8n+1,1}\cong U\oplus E_8^{\oplus n}$; see \cite{MilnorHusemoller1973}
and \cite[Chapter~V, \S2.2, Theorems~5--6]{Serre}.
We use the following unimodular gluing description.
\begin{lemma}[{\cite[Proposition~1.6.1]{Nikulin}}]
\label{gen:lem:nikulin}
Let $B$ be a primitive sublattice of an even unimodular
lattice $L$, and put $C=B^\perp\cap L$. Then there is an anti-isometry
$\phi:(A_B,q_B)\longrightarrow(A_C,q_C)$ such that
\[
 L=\{b+c\in B'\oplus C':\; [c]=\phi([b])\}.
\]
In particular, $|A_B|=|A_C|$.
\end{lemma}

For vectors of positive norm, we use the following orbit theorem.

\begin{lemma}[{\cite[Theorem~3.9.1]{Borcherds1984}}]
\label{gen:lem:positive-normal}
Let $n\ge1$, and let $y\in \II_{8n+1,1}$ be primitive with $y^2=2d>0$.
Then there is an isometry of $\II_{8n+1,1}$ that sends $y$ to $e-df$ in
$U\oplus E_8^{\oplus n}$.
\end{lemma}

We also need two scalar theta identities, both to construct roots
and to count decompositions in a hyperbolic summand. Let $P$ be an even
positive-definite unimodular lattice of rank $8$ or $16$, and put
\[
 \Theta_P(\tau)=\sum_{x\in P}q^{x^2/2},\qquad q=e^{2\pi i\tau}.
\]
The standard theta identities \cite[Chapter~VII, \S6.6]{Serre} give
\begin{equation}\label{gen:eq:small-theta}
 \Theta_P(\tau)=
 \begin{cases}
  E_4(\tau)=1+240\displaystyle\sum_{m\ge1}\sigma_3(m)q^m,
       &\rk P=8,\\[1mm]
  E_8(\tau)
       =1+480\displaystyle\sum_{m\ge1}\sigma_7(m)q^m,
       &\rk P=16,
 \end{cases}
\end{equation}
where $\sigma_j(m)=\sum_{a\mid m}a^j$, and $E_j(\tau)$ denotes the
normalized Eisenstein series of weight $j$. In particular, $E_8$ as
a lattice represents every positive even integer, and every $P$ of
either rank contains roots. These facts settle the vectors of nonnegative norm.

\begin{proposition}\label{gen:prop:nonnegative}
Let $n\ge1$ and $\gamma\in \II_{8n+1,1}$.
\begin{enumerate}
\item If $\gamma^2>0$ or $\gamma=0$, then
$R_n(\gamma)=\infty$.
\item Suppose $\gamma\ne0$ and $\gamma^2=0$. Write $\gamma=ky$,
where $k\ge1$ and $y$ is primitive. The quotient
\[
                    P_y=(y^\perp\cap \II_{8n+1,1})/\Z y
\]
is an even positive-definite unimodular of rank $8n$, and
\[
 R_n(\gamma)=
 \begin{cases}
   \infty,&P_y\text{ contains a root},\\
   0,&P_y\text{ contains no root}.
 \end{cases}
\]
\end{enumerate}
Consequently, for $n=1,2$, every vector of nonnegative norm has
infinitely many ordered decompositions into two roots.
\end{proposition}

\begin{proof}
Put $M=\II_{8n+1,1}$. If $\gamma^2>0$,
Lemma~\ref{gen:lem:positive-normal} reduces to $\gamma=k(e-df)$
with $k,d\ge1$. For each integer $a\ge k$,
Equation~\eqref{gen:eq:small-theta} provides $z_a\in E_8\subset E_8^{\oplus n}$ with
$z_a^2=2+2da(a-k)$, and
\begin{equation}\label{gen:eq:positive-pairs}
 r_a=ae+d(a-k)f+z_a,\qquad
 s_a=(k-a)e-da f-z_a
\end{equation}
give root decompositions of $\gamma$, distinct as $a$ varies since
the $e$-coefficient of $r_a$ is $a$. For $\gamma=0$, fix a root $z\in E_8$ and use
$(ae+z,-ae-z)$, with $a\in\Z$.

If $\gamma=ky\ne0$ is isotropic with $y$ primitive, unimodularity
gives $f_0\in M$ with $(y,f_0)=-1$; putting
$f=f_0+(f_0^2/2)y$ yields $f\in M$, $f^2=0$, and $(y,f)=-1$.
Thus $M=(\Z y+\Z f)\oplus P$, where $P\cong P_y$ is positive
definite even unimodular of rank $8n$. Any root decomposition
$r+s=ky$ satisfies $(r,y)=0$, so all such pairs are exactly
\begin{equation}\label{gen:eq:null-pairs}
 (r,s)=(ay+z,(k-a)y-z),\qquad
 a\in\Z,\quad z\in P,\quad z^2=2.
\end{equation}
This gives the stated criterion and infinitude. 
\end{proof}

\subsection{Gluing and primitive sublattices}

For negative norm, the perpendicular lattice is positive definite.
We express root decompositions as theta coefficients of this lattice
and then compare them with coefficients of a sublattice.

For a positive definite even lattice $L$, let $(\ee_\mu)_{\mu\in A_L}$
be the standard basis of $\C[A_L]$ and write
\[
 \Theta_L(\tau)=\sum_{\mu\in A_L}\sum_{v\in L+\mu}
 q^{Q(v)}\ee_\mu,\qquad q=e^{2\pi i\tau}.
\]
For a vector-valued Fourier series $F$, let $a_F(\mu,t)$ denote
the coefficient of $q^t\ee_\mu$. For a lattice, we abbreviate
\[
 a_L(\mu,t)=a_{\Theta_L}(\mu,t)
          =\#\{v\in L+\mu:Q(v)=t\}.
\]
Here $L+\mu$ is the corresponding coset in $L'$, and
$Q(\mu)=q_L(\mu)\in\Q/\Z$. A rational index $t$ is compatible if
$t-Q(\mu)\in\Z$; the lattice coefficient $a_L(\mu,t)$ is zero otherwise. For a unimodular
lattice, $A_L=\{0\}$ and we identify $\C[A_L]=\C\ee_0$ with $\C$.

\begin{lemma}\label{r9:lem:gluing-shell}
Let $M=\II_{8n+1,1}$ with $n\ge1$, and let $\gamma=ky\in M$, where $k\ge1$, $y\in M$ is primitive, and
\[
 y^2=-\delta<0,\qquad T=k^2\delta,
 \qquad K=y^\perp\cap M.
\]
Then $K$ is even positive definite of rank $8n+1$,
$\det K=\delta$, and $A_K$ is cyclic. There is a class $\mu\in A_K$
of order at most two such that
\begin{equation}\label{r9:eq:shell-bijection}
 R_n(\gamma)
 =\#\{v\in K+\mu:v^2=T/4+2\}
 =a_K(\mu,T/8+1).
\end{equation}
\end{lemma}
\begin{proof}
Apply Lemma~\ref{gen:lem:nikulin} to $\Z y\subset M$.
Its discriminant group is cyclic of order $\delta$, and the gluing
anti-isometry gives
\[
 \mu=\phi([ky/2]),\qquad 2\mu=0,
 \qquad q_K(\mu)=k^2\delta/8+\Z.
\]
The vector $ky/2$ belongs to $(\Z y)'$ because $\delta$ is even.
For $r+s=\gamma$ with $r^2=s^2=2$, subtraction gives
\[
 (r,\gamma)=\gamma^2/2=-T/2.
\]
Consequently $r=ky/2+v$ with $v\perp y$; gluing explains precisely
$v\in K+\mu$, and
\[
             r^2=-T/4+v^2=2.
\]
In contrast, such $v$ gives
$r=ky/2+v\in M$ and $s=\gamma-r=ky/2-v\in M$, both square two.
These constructions are inverse. Since $K$ is positive definite, $R_n(\gamma)$ is finite. 
\end{proof}

To estimate the coefficient in \eqref{r9:eq:shell-bijection}, we shall
use a fiber over a primitive sublattice of codimension one.
The following construction does not require a Lorentzian embedding.

\begin{lemma}\label{r9:lem:fibre}
Let $L$ be an even positive definite lattice of rank $r\ge2$ and
determinant $D$, and let $S\subset L$ be a primitive sublattice of
rank $r-1$ and determinant $c$. If $u$ is a primitive generator of
$L'\cap S_\R^\perp$, then $c=Du^2$. In particular,
\begin{equation}\label{r9:eq:n6}
 \lambda_1(L')^2=\frac1D
 \min_{\substack{S\subset L\ {\rm primitive}\\\rk S=r-1}}\det S,
\end{equation}
where $\lambda_1(L')$ is the length of a shortest nonzero vector of $L'$.
For every $\mu\in A_L$, there exist $\mu_S\in A_S$ and
$w\in S_\R^\perp$ such that
\begin{equation}\label{lattice:eq:general-fibre}
 w^2\le\frac D{4c},\qquad
 (L+\mu)\cap(w+S_\R)=w+(S+\mu_S).
\end{equation}
Consequently, for every compatible index $(\mu,t)$,
\[
 a_L(\mu,t)\ge a_S(\mu_S,t-Q(w)),\qquad
 t-Q(w)-Q(\mu_S)\in\Z.
\]
The sublattice $S$ does not need to have the least determinant, and the order of $\mu$ is unrestricted.
\end{lemma}
\begin{proof}
Since $u$ is primitive in $L'$, the functional
$(u,\cdot):L\to\Z$ is surjective. Choose $z\in L$
with $(u,z)=1$. Then $L=S+\Z z$, and the projection of $z$ onto
$S_\R^\perp$ is $w_0=u/u^2$. Taking determinants gives
\[
 D=\det S\cdot w_0^2=c/u^2.
\]
Every primitive corank-one sublattice arises in this way, and a
shortest nonzero dual vector is primitive. This proves
\eqref{r9:eq:n6}.

The projection of $L+\mu$ is given by $(\alpha+\Z)w_0$, with $\alpha\in \mathbb{Q}$. Its nearest
point $w$ to the origin satisfies $w^2\le w_0^2/4=D/(4c)$.
Choose $v_0\in L+\mu$ projected to $w$. Then
$x_0:=v_0-w\in S'$, since $(x_0,S)=(v_0,S)\subset\Z$.
Put $\mu_S:=x_0+S$. Primitivity gives $L\cap S_\R=S$, and hence
\[
 (L+\mu)\cap(w+S_\R)=v_0+S=w+(S+\mu_S).
\]
For $x\in S+\mu_S$, orthogonality gives $Q(w+x)=Q(w)+Q(x)$.
The coefficient inequality and compatibility follow.
\end{proof}

For the lattice $K$ in Lemma~\ref{r9:lem:gluing-shell}, we choose
$S\subset K$ of least determinant and write $c=\det S$. Then
$\lambda_1(K')^2=c/\delta$. Applying Lemma~\ref{r9:lem:fibre} to
the index $T/8+1$ gives
\begin{equation}\label{r9:eq:fibre-count}
 R_n(\gamma)\ge a_S(\mu_S,t),\qquad
 t=\frac T8+1-Q(w)\ge\frac T8\left(1-\frac1c\right)+1,
 \qquad t-Q(\mu_S)\in\Z,
\end{equation}
where we used $\delta\le T$. The fiber inequality holds for every
primitive $S$; the minimizing choice will be used to control
$\lambda_1(K')$ in Section~\ref{sec:lorentzian-estimates}.

\subsection{Unimodular summands and split vectors}

We first treat $c=1$. In this case a unimodular summand splits off,
and a scalar theta series gives the exact count.

\begin{proposition}\label{gen:prop:split-vectors}

Let $M=\II_{8n+1,1}$, let $\gamma=ky$ with $k\ge1$ and $y\in M$ primitive
of negative norm, and let $c$ be the least determinant of a primitive
rank-$8n$ sublattice of $y^\perp\cap M$. The following are equivalent:
\begin{enumerate}
\item $c=1$;
\item $y^\perp\cap M$ contains an even unimodular sublattice of rank $8n$;
\item there are a decomposition $M=U\oplus P$ and positive coprime
integers $a_0,b_0$ such that $P$ is even positive-definite unimodular and
$y=a_0e+b_0f$.
\end{enumerate}
In such a decomposition, the scalar theta coefficient is
\[
 a_P(0,m)=[q^m]\Theta_P=\#\{z\in P:Q(z)=m\}\qquad(m\in\Z_{\ge0}),
\]
and
\begin{equation}\label{gen:eq:split-root-count}
 R_n\bigl(k(a_0e+b_0f)\bigr)
 =\sum_{j=0}^k a_P\bigl(0,1+a_0b_0j(k-j)\bigr).
\end{equation}
For $k=1$, the first summands in all ordered decompositions are exactly
\begin{equation}\label{gen:eq:split-first-summands}
 \{a_0e+z:z\in P,\ z^2=2\}
 \ \sqcup\ \{b_0f+z:z\in P,\ z^2=2\}.
\end{equation}
\end{proposition}

\begin{proof}
A unimodular sublattice $P\subset M$ splits orthogonally:
the projection of $x\in M$ to $P_\Q$ belongs to $P'=P$. Thus
\[
 M=P\oplus(P^\perp\cap M),\qquad P^\perp\cap M\cong U.
\]
Since $y\perp P$, write $y=a_0e+b_0f$. Primitiveness gives
$\gcd(a_0,b_0)=1$, and $y^2=-2a_0b_0<0$ allows $a_0,b_0>0$ after changing
both signs. In contrast, the shown $P$ is primitive in
$y^\perp\cap M$ and has determinant one. This proves the equivalences.

For $r=pe+qf+z$, equations $r^2=(\gamma-r)^2=2$ are equivalent to
\[
 b_0p+a_0q=ka_0b_0,\qquad Q(z)=1+pq. 
\]
Coprimality gives
\[
 p=a_0j,\qquad q=b_0(k-j),\qquad
 Q(z)=1+a_0b_0j(k-j),\qquad j\in\Z.
\]
If $j<0$ or $j>k$, the last quantity is at most
$1-a_0b_0(k+1)<0$. Hence $0\le j\le k$, proving
\eqref{gen:eq:split-root-count}. The case $k=1$ gives
\eqref{gen:eq:split-first-summands}.
\end{proof}

\begin{corollary}\label{gen:cor:small-split-counts}
In the setting of Proposition~\ref{gen:prop:split-vectors},
suppose $n=1$ or $2$. Then
\begin{equation}\label{gen:eq:small-split-count}
 R_n\bigl(k(a_0e+b_0f)\bigr)
 =240n\sum_{j=0}^k\sigma_{4n-1}\bigl(1+a_0b_0j(k-j)\bigr).
\end{equation}
For $k=1$, this is $480n$. For $k\ge2$, it is at least
$240n(3+2^{4n-1})$, with equality if and only if $k=2$ and $a_0=b_0=1$.
Thus the primitive split counts are $480,960$, and the least
nonprimitive split counts are $2640,62880$, respectively.
\end{corollary}
\begin{proof}
Equation~\eqref{gen:eq:small-theta} gives the formula.
The two endpoints contribute $480n$. Every interior argument is an
integer $m\ge2$, and
$\sigma_{4n-1}(m)\ge1+m^{4n-1}\ge1+2^{4n-1}$. The equality holds only when there is exactly one interior term and its argument is $2$, which is equivalent to $k=2$ and $a_0b_0=1$. 
\end{proof}

\section{Theta coefficient estimates}\label{sec:theta-estimates}

Throughout this section, $L$ is an even positive definite lattice of
rank $r>4$. Put $\kappa=r/2>2$ and $D=\det L$.
The restriction $r>4$ ensures absolute convergence of the Eisenstein
and Poincar\'e series and of the Bessel sums after taking absolute
values. At weight two, the latter argument requires cancellation in
Kloosterman sums; see the discussion after
Lemma~\ref{r9:lem:poincare-bound}.
Additional hypotheses are stated in the individual results.
In particular, the uniform positive Eisenstein lower bound follows
under the local hypothesis in Lemma~\ref{r9:lem:local-lower}.
Section~\ref{sec:lorentzian-estimates} treats cyclic discriminant groups
and then applies the resulting bounds to Lorentzian perpendiculars.

\subsection{Eisenstein and cusp coefficients}\label{r9:sec:modular}

We write
\[
 e(z)=\exp(2\pi iz),\qquad
 \tau=x+iy\in\HH.
\]
The Hermitian form $(\ ,\ )_{\C[A_L]}$ on $\C[A_L]$ makes the standard basis
$(\ee_\mu)_{\mu\in A_L}$ orthonormal and is linear in the first
variable. We use the Weil representation $\rho_L$ determined by
\[
 \rho_L(\mathcal T)\ee_\mu=e(Q(\mu))\ee_\mu,\qquad
 \rho_L(\mathcal S)\ee_\mu=
 \frac{e(-r/8)}{\sqrt D}\sum_{\beta\in A_L}e(- (\mu,\beta))\ee_\beta.
\]
Here $\mathcal T$ and $\mathcal S$ are the standard metaplectic lifts of
$\left(\begin{smallmatrix}1&1\\0&1\end{smallmatrix}\right)$ and
$\left(\begin{smallmatrix}0&-1\\1&0\end{smallmatrix}\right)$, with
square-root factors $1$ and $\sqrt\tau$, respectively, and
$(\mu,\beta)=Q(\mu+\beta)-Q(\mu)-Q(\beta)\in\mathbb Q/\mathbb Z$.
For $\widetilde g=(g,\phi)\in\operatorname{Mp}_2(\mathbb Z)$, with
$g=\left(\begin{smallmatrix}a&b\\c&d\end{smallmatrix}\right)$ and
$\phi(\tau)^2=c\tau+d$, put
\[
 (F|_{\kappa,\rho_L}\widetilde g)(\tau)
 =\phi(\tau)^{-2\kappa}\rho_L(\widetilde g)^{-1}F(g\tau).
\]
Write $M_\kappa(\rho_L)$ and $S_\kappa(\rho_L)$ for the spaces
of holomorphic modular forms and cusp forms of weight $\kappa$,
respectively. Poisson summation gives
$\Theta_L\in M_\kappa(\rho_L)$, with constant term $\ee_0$.
We retain the coefficient notation $a_F(\mu,t)$ from
Section~\ref{sec:lattice-reductions}; its index is compatible when
$t-Q(\mu)\in\Z$.
All forms considered below satisfy $F_\mu=F_{-\mu}$, as follows from
the action of the center at weight $\kappa=r/2$.

Write $\widetilde\Gamma_\infty$ for the subgroup generated by $\mathcal T$ and
the center, and set
\[
 E_L=\sum_{\widetilde\Gamma_\infty\backslash\operatorname{Mp}_2(\mathbb Z)}
 \ee_0|_{\kappa,\rho_L}\widetilde g,
 \qquad F_L=\Theta_L-E_L.
\]
The series converges absolutely since $\kappa>2$, has constant term
$\ee_0$, and $F_L\in S_\kappa(\rho_L)$. Norms and pairings are taken
over
\[
 \mathcal F=\{\tau\in\HH:|x|\le\tfrac12,\ |\tau|\ge1\},\qquad
 \langle F,G\rangle=\int_{\mathcal F}(F(\tau),G(\tau))_{\C[A_L]}\,
 y^\kappa\frac{dx\,dy}{y^2}.
\]
The same integral defines $\|\Theta_L-\ee_0\|$, although the
function $\Theta_L-\ee_0$ is not modular.

We write $\Gamma$ and $\zeta$ for Euler's gamma function and
the Riemann zeta function, respectively.

\begin{lemma}\label{r9:lem:eisenstein-density}
For $\mu\in A_L$ and $t\in\Q_{>0}$ with $t-Q(\mu)\in\Z$,
\begin{equation}\label{r9:eq:eisenstein-product}
 a_{E_L}(\mu,t)=
 \frac{(2\pi)^\kappa t^{\kappa-1}}{\Gamma(\kappa)\sqrt D}
 \prod_p\alpha_p(L,\mu,t),
\end{equation}
where, for any representative $\widetilde\mu\in L'$,
\begin{align*}
 N_{p^a}(L,\mu,t)
 &=\#\{v\in L/p^aL:Q(v+\widetilde\mu)-t\equiv0\pmod{p^a}\},\\
 \alpha_p(L,\mu,t)
 &=\lim_{a\to\infty}p^{a(1-r)}N_{p^a}(L,\mu,t).
\end{align*}
The local limits exist, and the product is absolutely convergent.
\end{lemma}

\begin{proof}
Use the coefficient formula of Bruinier and Kuss \cite{BK} in the
normalization of \cite[Proposition~3.2.1 and (3.2.7)--(3.2.9)]{Opitz},
whose constant term is $\ee_0$ rather than $2\ee_0$.
Apply it to $L(-1)$ at weight $\kappa=r/2$; the dual Weil
representation of $L(-1)$ is our $\rho_L$. The archimedean sign is
$(-1)^{(2\kappa-r)/4}=1$.
After $v\mapsto-v$, the congruence counts agree with
$N_{p^a}(L,\mu,t)$; \cite[Lemma~3.2.2]{Opitz} identifies their
stabilized normalized value with the local factor at $s=\kappa$.
This proves \eqref{r9:eq:eisenstein-product}.
At all but finitely many primes,
$\alpha_p=1+O_r(p^{1-r/2})$ by \cite[Theorem~3.2.3]{Opitz},
so the product converges absolutely.
\end{proof}

\begin{lemma}\label{r9:lem:local-lower}
Let $p$ be a prime. Suppose that $L\otimes\Z_p$ contains a
unimodular orthogonal summand of rank $2m$, with $m\in\Z_{\ge1}$, which
is even if $p=2$. For every compatible positive index $(\mu,t)$,
\begin{equation}\label{gen:eq:local-density-lower}
 \alpha_p(L,\mu,t)\ge
 1-(p-1)p^{-m}-p^{1-2m}.
\end{equation}
In particular,
\[
 \alpha_p(L,\mu,t)\ge1-p^{1-m}\quad(p\ne2),\qquad
 \alpha_2(L,\mu,t)\ge1-2^{-m}-2^{1-2m}.
\]
If the hypothesis holds at every prime with the same $m\ge3$, then
\begin{equation}\label{gen:eq:eta-m}
 \prod_p\alpha_p(L,\mu,t)\ge
 \eta(m):=\frac{1+2^{-m}}{\zeta(m-1)}.
\end{equation}
The parameters used below are
\[
 \eta(3)=\frac9{8\zeta(2)},\qquad
 \eta(4)=\frac{17}{16\zeta(3)},\qquad
 \eta(m)>\frac{49}{50}\quad(m\ge7).
\]
\end{lemma}

\begin{proof}
Write $L_p=U_p\oplus C_p$, with $U_p$ unimodular of rank $2m$.
Write $v=v_U+v_C$, with $v_U\in U_p$ and $v_C\in C_p$, and decompose
\[
 \widetilde\mu=\widetilde\mu_U+\widetilde\mu_C
 \in U_p'\oplus C_p'.
\]
Since $U_p'=U_p$, the change of variables
\[
 u:=v_U+\widetilde\mu_U
\]
permutes the residue classes in $U_p/p^aU_p$.
For each fixed $v_C\pmod{p^aC_p}$, orthogonality rewrites the
congruence $Q(v+\widetilde\mu)\equiv t\pmod{p^a}$ defining
$N_{p^a}(L,\mu,t)$ as
\[
 Q(u)\equiv t_0\pmod{p^a},\qquad
 t_0:=t-Q(v_C+\widetilde\mu_C).
\]
Moreover,
\[
 t_0=\bigl(t-Q(\widetilde\mu)\bigr)
     +Q(\widetilde\mu_U)-Q(v_C)
     -(v_C,\widetilde\mu_C)\in\Z_p.
\]

The quadratic form on $U_p/pU_p$ is nondegenerate, also at $p=2$.
For some $\varepsilon_p\in\{\pm1\}$ its representation numbers
\cite[\S3, (4)--(5)]{Casselman2018} are
\[
 N_p(t_0)=
 \begin{cases}
 p^{2m-1}+\varepsilon_p(p-1)p^{m-1},&t_0\equiv0\pmod p,\\
 p^{2m-1}-\varepsilon_p p^{m-1},&t_0\not\equiv0\pmod p.
 \end{cases}
\]
These formulas give at least
$p^{2m-1}-(p-1)p^{m-1}-1$ nonzero solutions in every fiber. 
Every nonzero vector is nonsingular, and each solution has exactly
$p^{(a-1)(2m-1)}$ lifts modulo $p^a$;
see \cite[Lemma~3.2]{Hanke2004}.
Counting the $p^{a(r-2m)}$ choices of the complementary coordinate
$v_C$ therefore gives
\[
 p^{a(1-r)}N_{p^a}(L,\mu,t)
 \ge p^{1-2m}\bigl(p^{2m-1}-(p-1)p^{m-1}-1\bigr)
 =1-(p-1)p^{-m}-p^{1-2m}.
\]
This proves \eqref{gen:eq:local-density-lower}.

For odd $p$, $p^{1-2m}\le p^{-m}$ gives the stated weaker bound.
At $p=2$, use
\[
 1-2^{-m}-2^{1-2m}=(1+2^{-m})(1-2^{1-m}).
\]
Multiplication over primes proves \eqref{gen:eq:eta-m}.
The factors increase with $m$, and $\zeta(6)<327/320$, so
\[
 \eta(m)\ge\eta(7)>
 \frac{129}{128}\frac{320}{327}=\frac{215}{218}>\frac{49}{50}
 \qquad(m\ge7).
\]
\end{proof}

\begin{remark}\label{r9:rem:local-hypotheses}
For an even lattice of rank $r$, the rank of a Gram matrix modulo $p$
is $r-\ell(A_L/pA_L)$. Its unimodular Jordan component has this
rank; at $p=2$ it is even. If $S\subset K$ is primitive of codimension
one and $\ell(A_K)\le\ell$, a basis of $S$ extends to one of $K$.
Deleting the last row and column decreases the rank modulo $p$ by
at most two, giving $\ell(A_S)\le\ell+1$.
\end{remark}

We next turn to the cuspidal part. The following estimate reduces
its norm to a theta norm that can be bounded by lattice-point counts.

\begin{lemma}\label{r9:lem:cusp-projection}
Define
\[
 c_\kappa^2=
 \frac2{\kappa-1}\int_0^{1/2}
       \left[1-(1-x^2)^{(\kappa-1)/2}\right]dx.
\]
Then
\[
 \|F_L\|\le\|\Theta_L-\ee_0\|+c_\kappa,
 \qquad
 c_\kappa\le
 \frac{(\sqrt3/2)^{\min\{\kappa-3,0\}/2}}{\sqrt{24}}.
\]
In particular, $c_\kappa\le1/\sqrt{24}$ for every $r\ge6$.
\end{lemma}
\begin{proof}
Since $\kappa>2$, the Petersson pairing is absolutely convergent,
and unfolding gives $\langle E_L,F_L\rangle=0$. The component $F_{L,0}$ has zero average over a full period.
Put $g=\ee_0\mathbf1_{y<1}$. Since $\Theta_L-\ee_0$ and $g$ belong to
$L^2(\mathcal F,y^{\kappa-2}dx\,dy;\C[A_L])$,
the Cauchy--Schwarz inequality gives
\[
 \|F_L\|^2
 =\langle\Theta_L-\ee_0,F_L\rangle+\langle g,F_L\rangle
 \le\bigl(\|\Theta_L-\ee_0\|+\|g\|\bigr)\|F_L\|,
\]
where
\[
 \|g\|^2=
 \int_{-1/2}^{1/2}\int_{\sqrt{1-x^2}}^1y^{\kappa-2}\,dy\,dx
 =c_\kappa^2.
\]
Since $\sqrt3/2\le y\le1$ on the support of $g$,
\[
 c_\kappa^2=\|g\|^2
 =\int_{-1/2}^{1/2}\int_{\sqrt{1-x^2}}^1
       y^{\kappa-2}\,dy\,dx
\le\frac{(\sqrt3/2)^{\min\{\kappa-3,0\}}}{2}
       \int_{-1/2}^{1/2}x^2\,dx
=\frac{(\sqrt3/2)^{\min\{\kappa-3,0\}}}{24}.
\]
Dividing by $\|F_L\|$ proves the assertion when $F_L\ne0$;
the zero case is immediate.
\end{proof}

To estimate Fourier coefficients from the cusp norm, we retain the
dependence on the discriminant group in the Poincar\'e-series formula.
For a finite abelian group $A$ and an integer $a$, put
$A[a]=\{x\in A:ax=0\}$. Let $N_A$ be its exponent, with $N_A=1$
for the trivial group, and define
\[
 \mathcal H(A)=\frac1{N_A}\sum_{a=1}^{N_A}\sqrt{|A[a]|}.
\]

\begin{lemma}\label{r9:lem:weil-bound}
For $g=\left(\begin{smallmatrix}a&b\\c&d\end{smallmatrix}\right)\in\SL_2(\Z)$,
any lift $\widetilde g$, and $\mu,\beta\in A_L$,
\[
 |\rho_L(\widetilde g)_{\mu,\beta}|
 \le\sqrt{\frac{|A_L[c]|}{D}}.
\]
\end{lemma}

\begin{proof}
By Str\"omberg's explicit formula
\cite[Theorem~6.4 and Remark~6.8]{Stromberg2013}, each matrix entry
has modulus either $0$ or $\sqrt{|A_L[c]|/D}$.
Changing the metaplectic lift multiplies the matrix by a scalar of
modulus one.
\end{proof}

\begin{lemma}\label{r9:lem:poincare-bound}
Put $\nu=\kappa-1>1$.
For every $F\in S_\kappa(\rho_L)$ and every compatible positive
index $(\mu,t)$, one has
\begin{equation}\label{r9:eq:general-cusp-monotone}
 |a_F(\mu,t)|\le
 \frac{\|F\|(4\pi t)^{\nu/2}}{\sqrt{\Gamma(\nu)}}
 \left(1+C(\kappa)\frac{\mathcal H(A_L)t}{\sqrt D}\right)^{1/2},
\end{equation}
where
\[
 M_\nu=\frac{27}{40\nu^{1/3}},\qquad
 C(\kappa)=
 \frac{4\pi^2\nu}{(\nu-1)\Gamma(\nu+1)^{1/\nu}}
              M_\nu^{1-1/\nu}.
\]
In particular, the following upper bounds hold in the four ranks
used below:
\begin{equation}\label{r17:eq:cuspcoef}
 \begin{array}{c|cccc}
 r&8&9&16&17\\\hline
 \text{upper bound for }C(\kappa)&20&16&6&6
 \end{array}.
\end{equation}
\end{lemma}

\begin{proof}
Put $f(c)=\sqrt{|A_L[c]|}$, and let $\mathcal P_{\mu,t}$ be the
Poincar\'e series with seed $e(t\tau)(\ee_\mu+\ee_{-\mu})/2$.
Unfolding the Petersson pairing as in
\cite[\S1.2]{Bruinier2002} gives, in the normalization of
\cite[\S3, (8) and Lemma~3.1]{BBFW2026},
\[
 \langle F,\mathcal P_{\mu,t}\rangle
 =\frac{\Gamma(\nu)}{(4\pi t)^\nu}a_F(\mu,t),\qquad
 \|\mathcal P_{\mu,t}\|^2
 =\frac{\Gamma(\nu)}{(4\pi t)^\nu}
       a_{\mathcal P_{\mu,t}}(\mu,t).
\]
In particular, $a_{\mathcal P_{\mu,t}}(\mu,t)\ge0$.
The Cauchy--Schwarz inequality therefore gives
\begin{equation}\label{r9:eq:poincare-cauchy}
 |a_F(\mu,t)|
 =\frac{(4\pi t)^\nu}{\Gamma(\nu)}
       |\langle F,\mathcal P_{\mu,t}\rangle|
 \le\frac{(4\pi t)^\nu}{\Gamma(\nu)}
       \|F\|\,\|\mathcal P_{\mu,t}\|
 =\frac{\|F\|(4\pi t)^{\nu/2}}{\sqrt{\Gamma(\nu)}}
       \sqrt{a_{\mathcal P_{\mu,t}}(\mu,t)}.
\end{equation}

In the Kloosterman--Bessel expansion
\cite[proof of Lemma~3.2]{BBFW2026}, each Kloosterman sum has at
most $c$ terms, each of modulus at most $f(c)/\sqrt D$ by
Lemma~\ref{r9:lem:weil-bound}.
Since the diagonal seed term is at most $1$, we obtain
\begin{equation}\label{r9:eq:poincare-bessel}
 0\le a_{\mathcal P_{\mu,t}}(\mu,t)
 \le1+\frac{2\pi}{\sqrt D}
       \sum_{c\ge1}f(c)|J_\nu(4\pi t/c)|,
\end{equation}
where $J_\nu$ is the Bessel function of the first kind.
Moreover, Landau's inequality $|J_\nu(x)|<M_\nu$
\cite{Landau2000} and the elementary power bound
\cite[(10.14.4)]{DLMF} imply
\[
 |J_\nu(x)|\le b_\nu(x):=
 \min\{M_\nu,(x/2)^\nu/\Gamma(\nu+1)\}.
\]
Observe that the function $b_\nu$ is nondecreasing. Thus, for $X>0$,
\[
 \sum_{m\ge1}b_\nu(X/m)
 \le\int_0^\infty b_\nu(X/t)\,dt
 =XI_\nu,\qquad
 I_\nu=\int_0^\infty b_\nu(x)x^{-2}\,dx
 =\frac{\nu M_\nu^{1-1/\nu}}
        {2(\nu-1)\Gamma(\nu+1)^{1/\nu}}.
\]
Let $h$ be multiplicative, with $h(1)=1$ and
\[
 h(p^j)=f(p^j)-f(p^{j-1})\ge0\qquad(j\ge1).
\]
It is supported on divisors of $N_{A_L}$, and
$f(c)=\sum_{a\mid c}h(a)$. Averaging over $1\le c\le N_{A_L}$ gives
$\mathcal H(A_L)=\sum_{a\mid N_{A_L}}h(a)/a$. Consequently, for the sum in \eqref{r9:eq:poincare-bessel} we get
\[
 \sum_{c\ge1}f(c)|J_\nu(X/c)|
 \le\sum_{a\ge1}h(a)\sum_{m\ge1}b_\nu(X/(am))
 \le XI_\nu\mathcal H(A_L).
\]
Substitute $X=4\pi t$ in
\eqref{r9:eq:poincare-bessel}, and note that
$8\pi^2I_\nu=C(\kappa)$. This completes the proof.
\end{proof}

The restriction $\kappa>2$ is essential for the absolute-value
estimate \eqref{r9:eq:poincare-bessel}: since
$J_\nu(X/c)\sim (X/2c)^\nu/\Gamma(\nu+1)$ as $c\to\infty$
and $\sqrt{|A_L[c]|}\ge1$, its series diverges for $0<\nu\le1$.
A weight-two estimate by this method therefore requires cancellation
in the Kloosterman sums.

The next estimates for $\mathcal H(A)$ will be used for the cyclic
group $A_K$ and the group $A_S$, which has at most two generators.

\begin{lemma}\label{r9:lem:H-bound}
For a finite abelian group $A$ of order $D$, write
$A_p=\bigoplus_i\Z/p^{e_{p,i}}\Z$. For any $j\geq 0$, define 
$t_{p,j}=\sum_i\min\{j,e_{p,i}\}$. Then
\[
 \mathcal H(A)=\prod_{p\mid D}
 \left(1+\sum_{j=1}^{\max_i e_{p,i}}
       \frac{p^{t_{p,j}/2}-p^{t_{p,j-1}/2}}{p^j}\right).
\]
In particular,
\begin{equation}\label{r17:eq:Hbounds}
 \mathcal H(A)<\tfrac98D^{1/6}\quad\text{if $A$ is cyclic},\qquad
 \mathcal H(A)<\tfrac98D^{1/4}\quad\text{if $\ell(A)\le2$}.
\end{equation}
\end{lemma}

\begin{proof}
Put $n_p=\max_i e_{p,i}$. Since $N_A=\prod_{p\mid D}p^{n_p}$,
the Chinese remainder theorem gives
\[
 \mathcal H(A)=\prod_{p\mid D}H_p,\qquad
 H_p=\frac1{p^{n_p}}\sum_{k\bmod p^{n_p}}\sqrt{|A_p[k]|}.
\]
Note that $|A_p[k]|=p^{t_{p,j}}$ if $p^j\mid k$ and $p^{j+1}\nmid k$. Summation by parts yields
\[
 H_p=1+\sum_{j=1}^{n_p}
 \frac{p^{t_{p,j}/2}-p^{t_{p,j-1}/2}}{p^j},
\]
proving the product formula.

If $\ell(A)\le2$, write
$A_p=\Z/p^e\Z\oplus\Z/p^f\Z$, with $0\le e\le f$.
The corresponding local factor is
\[
 H_p(e,f)=1+e(1-p^{-1})
             +p^{-1/2}(1-p^{-(f-e)/2}).
\]
The ratios $H_p(e+1,f+1)/H_p(e,f)$ and
$H_p(e,f+1)/H_p(e,f)$ decrease as $e$ or $f-e$
increases.

For cyclic groups, $e=0$. If $p\ge5$, then
\[
 \frac{H_p(0,1)}{H_p(0,0)}
 =1+p^{-1/2}-p^{-1}\le\frac54<p^{1/6},
\]
so $H_p(0,f)p^{-f/6}$ has maximum $1$.
For $p=2,3$,
\[
 \frac{H_2(0,2)}{H_2(0,1)}
 =\frac3{\sqrt2}-1<2^{1/6},\qquad
 \frac{H_3(0,2)}{H_3(0,1)}
 =4-\frac5{\sqrt3}<3^{1/6}.
\]
Thus only $f=0,1$ need be compared, and the maximum occurs
at $f=1$. Hence
\[
 \frac{\mathcal H(A)}{D^{1/6}}
 \le
 \frac{(1/2+1/\sqrt2)(2/3+1/\sqrt3)}{6^{1/6}}
 <\frac98.
\]

For $\ell(A)\le2$, normalize by $p^{(e+f)/4}$.
The two adjacent ratios must now be compared with $\sqrt p$
and $p^{1/4}$, respectively. For $p\ge3$, their initial values
satisfy
\[
 2-p^{-1}\le\sqrt p,\qquad
 1+p^{-1/2}-p^{-1}\le\frac54<p^{1/4},
\]
so the normalized maximum is $1$.
At $p=2$, the inequalities
\[
 \frac{H_2(2,2)}{H_2(1,1)}=\frac43<\sqrt2,\qquad
 \frac{H_2(0,2)}{H_2(0,1)}<2^{1/4}
\]
reduce the maximum to $e,f-e\in\{0,1\}$.
Comparison of these four values gives maximum
$3/(2\sqrt2)$ at $e=f=1$. Therefore
\[
 \frac{\mathcal H(A)}{D^{1/4}}
 \le\frac3{2\sqrt2}<\frac98.
\]
\end{proof}

\subsection{Norm estimates and theta coefficient bounds}
\label{r9:sec:norms}

Uniform Petersson-norm estimates for scalar cuspidal theta series are
given in \cite{Rouse2019,Waibel2022}. Here we estimate the full vector-valued
theta series after subtracting its constant term, with explicit
constants depending only on the rank and complementary determinants.

For a Euclidean lattice $B$ of rank $r$, define
\[
 \lambda_j(B)=\min\bigl\{R>0:\;
 \dim_\R\Span_\R\{u\in B:u^2\le R^2\}\ge j\bigr\}
 \qquad(1\le j\le r).
\]
Thus, $\lambda_j(B)$ is the smallest radius of a ball centered at the
origin containing $j$ linearly independent vectors of $B$. In
particular, $\lambda_1(B)^2=\min_{0\ne u\in B}u^2$. To bound $\|\Theta_L-\ee_0\|$, we count short vectors of $L'$
using the successive minima $\lambda_j(L')$ and estimate the
remaining tail separately. For convenience, we set
\begin{equation}
 C_r=\frac{\Gamma(r/2-1)}{(2\pi)^{r/2-1}},  
\end{equation}
and define
\begin{equation}\label{r9:eq:n1}
 I_r(N)=\int_{\sqrt3/2}^\infty y^{r/2-2}e^{-2\pi Ny}\,dy\qquad(N>0).
\end{equation}

For each theta series, the Petersson norm is taken at half the rank
of its lattice.

\begin{lemma}\label{r9:lem:det-bounds}
For $j\ge0$, define $\beta_j$ by periodicity modulo $8$ and
\begin{equation}\label{r9:eq:n3}
 (\beta_0,\ldots,\beta_7)=(1,2,3,4,4,4,3,2).
\end{equation}
Then every even positive definite lattice $J$ of rank $j$ satisfies
$\det J\ge\beta_j$. If $8\mid j$ and $\det J>1$, then $\det J\ge4$.
\end{lemma}
\begin{proof}
Milgram's formula \cite[Appendix~4, p.~127]{MilnorHusemoller1973} yields
\begin{equation}\label{r9:eq:n4}
 |A_J|^{-1/2}\sum_{\mu\in A_J}e^{2\pi iQ(\mu)}=e^{2\pi ij/8}.
\end{equation}
For discriminant orders $1,2,3$, the possible normalized Gauss sums are,
respectively, $1$, $e^{\pm2\pi i/8}$ and $e^{\pm2\pi i/4}$:
the quadratic forms are $0$, $Q(1)=\pm1/4$ on $\Z/2\Z$, and
$Q(1)=\pm1/3$ on $\Z/3\Z$. Thus these determinant orders require
$j\equiv0$, $j\equiv\pm1$ and $j\equiv\pm2\pmod8$, respectively.
These restrictions give the assertions, with
the convention $\det(0)=1$ in rank zero.
\end{proof}

\begin{lemma}\label{r9:lem:minima-products}
Let $L$ be an even positive-definite lattice of rank $r$ with $D=\det L$. For any $1\le j\le r$, we have
\begin{equation}\label{r9:eq:n5}
 \prod_{i=1}^j\lambda_i(L')\geq \sqrt{\beta_{r-j}/D}=:p_j. 
\end{equation}
\end{lemma}
\begin{proof}
We choose independent $u_i\in L'$ with $i=1,...,j$ such that $u_i^2=\lambda_i(L')^2$, and put $\Lambda=\sum_{i=1}^j\Z u_i$, $\overline\Lambda=L'\cap(\Lambda\otimes\Q)$ and
$J=L\cap\Lambda^\perp$. Note that $\det J \geq \beta_{r-j}$. Extend a basis of $\overline\Lambda$ to one of $L'$, and then take the dual basis of $L$,  whose last $r-j$ vectors form a basis of $J$. With respect to these bases, we compute the Gram matrices of $L'$, $J$, and $\overline\Lambda$, and thus conclude that   $\det\overline\Lambda=\det J/D$. Therefore, 
\[
 \prod_{i=1}^j\lambda_i(L')^2\ge\det\Lambda
 =[\overline\Lambda:\Lambda]^2\det\overline\Lambda
 \ge\frac{\beta_{r-j}}D.
\]
\end{proof}

We convert these product bounds into lattice-point bounds using
the ellipsoid theorem of Malikiosis \cite[Theorem~1.1]{Malikiosis2013}.
For a Euclidean lattice $B$ of rank $j$ and $R>0$, it gives
\begin{equation}\label{norm:eq:ball-count}
 \#\{v\in B:v^2\le R^2\}
 \le\prod_{i=1}^j\left\lfloor\frac{2R}{\lambda_i(B)}+1\right\rfloor.
\end{equation}

We also need uniform bounds for individual shells. For $j\ge0$, put
\[
 M_r(j)=\binom{r+j}{j}-\binom{r+j-2}{j-2},\qquad
 \mathcal U_r=1+2\sum_{m\ge1}\binom{r+4m-2}{4m-1}e^{-2\pi m},
\]
where the second binomial coefficient is zero for $j<2$. Define
\begin{equation}\label{norm:eq:general-tail-constant}
 \tau_r=\mathcal U_r\max_{N\in\{1/4\}\cup\{j/2:1\le j\le r\}}
       M_r(\lfloor2N\rfloor)e^{\pi N}I_r(N),
\end{equation}
where $I_r(N)$ is defined in \eqref{r9:eq:n1}.

\begin{lemma}\label{norm:lem:general}
Let $2\le j\le r$ and $\alpha\ge\tfrac12+\tfrac{r-2}{2j}$. Set
\begin{equation}\label{r9:eq:n19}
 C_{r,j,\alpha}=C_r3^j\beta_{r-j}^{-\alpha}
                         4^{(r-2)/2-j\alpha}.
\end{equation}
For any choices $\alpha_j\ge\tfrac12+\tfrac{r-2}{2j}$,
\begin{equation}\label{r9:eq:n20}
 \|\Theta_L-\ee_0\|^2
 \le2\zeta(r-2)C_r\lambda_1(L')^{2-r}
      +\sum_{j=2}^rC_{r,j,\alpha_j}D^{\alpha_j}+\tau_r\sqrt D.
\end{equation}
Consequently,
\begin{equation}\label{norm:eq:det-only}
 \|\Theta_L-\ee_0\|^2
 \le\mathcal B_r D^{(r-2)/2}+\tau_r\sqrt D,
 \quad
 \mathcal B_r=2\zeta(r-2)C_r\beta_{r-1}^{-(r-2)/2}
             +\sum_{j=2}^rC_{r,j,(r-2)/2}.
\end{equation}
For the four ranks needed below, one may use
\begin{equation}\label{r17:eq:tail}
 \tau_8,\tau_9<3,\qquad \tau_{16},\tau_{17}<2000.
\end{equation}
\end{lemma}
\begin{proof}
Since $\mathcal F$ is contained in the strip $|x|\le1/2$, $y\ge\sqrt{3}/2$,
Fourier orthogonality within each discriminant component gives
\begin{equation}\label{r9:eq:n24}
 \|\Theta_L-\ee_0\|^2
 \le\sum_{\mu\in A_L}
       \sum_{\substack{t>0\\t-Q(\mu)\in\Z}}
       a_L(\mu,t)^2 I_r(2t)
 =\sum_{0\ne v\in L'}a_L([v],Q(v))I_r(v^2).
\end{equation}
In the last equality, a shell of multiplicity $a_L(\mu,t)$ is
counted by its vectors, each carrying the same weight $a_L(\mu,t)$.
We split this sum into $0<v^2<1/4$ and $v^2\ge 1/4$. That is to say,
\begin{equation}\label{e_split_sum}
 \sum_{0\ne v\in L'}a_L([v],Q(v))I_r(v^2)=\sum_{\substack{v\in L'\\0<v^2<1/4}}a_L([v],Q(v))I_r(v^2)+\sum_{\substack{v\in L'\\v^2\ge 1/4}}a_L([v],Q(v))I_r(v^2). 
\end{equation}

First, we bound the second sum, i.e., the tail, in \eqref{e_split_sum}. For distinct vectors $v,w$ in the shell $v^2=w^2=N$ of a coset $L+\mu$,
\[
 N-(v,w)=(v-w)^2/2\in\{1,\ldots,\lfloor2N\rfloor\}.
\]
The absolute bound \cite[Theorem~4.8]{DGS1977} therefore gives
\begin{equation}\label{r9:eq:n8}
 a_L(\mu,N/2)\le M_r(\lfloor2N\rfloor),\qquad
 a_L(\mu,N/2)\le1\quad(N<1/2).
\end{equation}
For integers $m\ge1$, the integral shell $v^2=2m$ is antipodal,
and its normalized inner products between distinct vectors lie in
$\{j/(2m):-2m\le j\le2m-1\}$.
The antipodal absolute bound \cite[Theorem~6.8]{DGS1977} gives
\begin{equation}\label{r9:eq:n9}
 \#\{v\in L:v^2=2m\}\le2\binom{r+4m-2}{4m-1},
 \qquad\sum_{v\in L}e^{-\pi v^2}\le\mathcal U_r.
\end{equation}

Put $G_r(N)=e^{\pi N}I_r(N)$. Its defining integral gives
$(\log G_r)'\le-\pi(\sqrt3-1)$, and
\[
 \frac{M_r(j+1)}{M_r(j)}
 =\frac{r+2j+1}{r+2j-1}\frac{r+j-1}{j+1}
 <2<e^{\pi(\sqrt3-1)/2}\qquad(j\ge r).
\]
Thus $M_r(\lfloor2N\rfloor)G_r(N)$ decreases between consecutive
half-integers, and its values at $j/2$ decrease for $j\ge r$.
By \eqref{norm:eq:general-tail-constant},
\[
 M_r(\lfloor2N\rfloor)G_r(N)\le\frac{\tau_r}{\mathcal U_r}
 \qquad(N\ge 1/4).
\]
Combining this with \eqref{r9:eq:n8} and Poisson summation yields
\begin{equation}\label{r9:eq:n10}
 \sum_{\substack{v\in L'\\v^2\ge 1/4}}
       a_L([v],Q(v))I_r(v^2)
 \le\frac{\tau_r}{\mathcal U_r}
       \sum_{v\in L'}e^{-\pi v^2}
 =\frac{\tau_r\sqrt D}{\mathcal U_r}
       \sum_{v\in L}e^{-\pi v^2}
 \le\tau_r\sqrt D.
\end{equation}
Next, we estimate the first sum in \eqref{e_split_sum}. Put
\[
 B(R)=\#\{v\in L':0<v^2\le R^2\},\qquad R_j=p_j^{1/j}, \quad p_j=\sqrt{\beta_{r-j}/D}.
\]
If the vectors counted by $B(R)$ span a space $V$ of dimension $j\ge2$,
the first $j$ minima of $L'\cap V$ agree with those of $L'$ and
are at most $R$. Hence \eqref{norm:eq:ball-count} and
\eqref{r9:eq:n5} give
\[
 B(R)+1
 \le\prod_{i=1}^j
       \left\lfloor\frac{2R}{\lambda_i(L')}+1\right\rfloor
 \le\frac{3^jR^j}{p_j}.
\]
Moreover,
\[
 p_j\le\prod_{i=1}^j\lambda_i(L')\le R^j,
\]
so $R\ge R_j$.
In dimension one, $B(R)=2\lfloor R/\lambda_1(L')\rfloor$.
Taking an upper bound over all possible dimensions gives
\begin{equation}\label{r9:eq:n23}
 B(R)\le2\lfloor R/\lambda_1(L')\rfloor+
       \sum_{j=2}^r\frac{3^jR^j}{p_j}\mathbf1_{\{R\ge R_j\}}.
\end{equation}
For $0<v^2<1/4$, the multiplicity in \eqref{r9:eq:n24} is one
by \eqref{r9:eq:n8}. Stieltjes integration by parts gives
\[
 \sum_{\substack{v\in L'\\0<v^2<1/4}}I_r(v^2)
 \le I_r(1/4)B(1/2)
 +\int_0^{1/2}B(R)
       \left(-\frac{d}{dR}I_r(R^2)\right)\,dR.
\]
Extending the defining positive integrals to $y\ge0$ gives
\[
 I_r(R^2)\le C_rR^{2-r},\qquad
 -\frac{d}{dR}I_r(R^2)\le(r-2)C_rR^{1-r}.
\]
Substitute \eqref{r9:eq:n23} in the preceding integral bound.
The first term of the counting bound contributes
\[
 2\sum_{\substack{m\ge1\\m\lambda_1(L')\le 1/2}}
       I_r\bigl((m\lambda_1(L'))^2\bigr)
 \le2\zeta(r-2)C_r\lambda_1(L')^{2-r}.
\]

For $j\ge2$, the corresponding term vanishes if $R_j>1/2$.
Otherwise, it is at most
\begin{equation}\label{e_first_sum_others}
  \frac{C_r3^j}{p_j}
 \left(2^{r-j-2}+(r-2)\int_{R_j}^{1/2}R^{j-r+1}\,dR\right).
\end{equation}
Set
\[
 u=\frac{D}{4^j\beta_{r-j}}
   =\frac1{(2R_j)^{2j}}\ge1,
 \qquad b=\frac{r-2}{2j}.
\]
The substitution $t=(2R)^{-2j}$ rewrites this bound \eqref{e_first_sum_others} as
\[
 C_{r,j,\alpha}D^\alpha u^{1/2-\alpha}
 \left(1+b\int_1^u t^{b-3/2}\,dt\right).
\]
For $t\ge1$, we have $t^{b-3/2}\le t^{b-1}$, so
\[
 u^{1/2-\alpha}
 \left(1+b\int_1^u t^{b-3/2}\,dt\right)
 \le u^{1/2+b-\alpha}\le1
 \qquad\left(\alpha\ge\frac12+b\right).
\]
Applying this with $\alpha=\alpha_j$ and adding the tail estimate
\eqref{r9:eq:n10} proves \eqref{r9:eq:n20}.
Taking $\alpha_j=(r-2)/2$ and using
$\lambda_1(L')^2\ge\beta_{r-1}/D$ gives \eqref{norm:eq:det-only}.

Finally, for the four ranks in \eqref{r17:eq:tail}, put
$q_r=\binom{r+6}{7}/\binom{r+2}{3}$. The ratios of consecutive
binomial coefficients in $\mathcal U_r$ decrease, so
\[
 \mathcal U_r\le
 1+\frac{2\binom{r+2}{3}e^{-2\pi}}{1-q_re^{-2\pi}}.
\]
For $h=r/2-2$, $b=\lceil h\rceil$ and $x=2\pi N$, integration by
parts gives
\begin{equation}\label{r9:eq:n16}
 G_r(N)\le\frac{e^{-(\sqrt3-1)\pi N}}x
        \sum_{\ell=0}^b\frac{(h)_\ell (\sqrt3/2)^{h-\ell}}{x^\ell}.
\end{equation}
Here $(h)_\ell=h(h-1)\cdots(h-\ell+1)$ and $(h)_0=1$.
The remainder is nonpositive because $(h)_{b+1}\le0$.
Substitution at the finitely many indices in
\eqref{norm:eq:general-tail-constant}, evaluated using rational intervals,
gives \eqref{r17:eq:tail}.
\end{proof}

\begin{corollary}\label{norm:cor:thick-thin}
Let $r\ge8$, $T>0$ and $b>0$. Let $\mathcal B_{r}$ be defined in \eqref{norm:eq:det-only}. Put
\[
 \mathcal D_r(b)=C_r b^{-(r-2)/2}
       \left(\frac{4(r-2)}{r-4}+\frac{4(r-2)}{r-3}\right)
       +\sum_{j=3}^rC_{r,j,(r-2)/4}.
\]
\begin{enumerate}
\item If $K$ is even positive definite of rank $r$, with
$\det K\le T$ and $\lambda_1(K')^2\ge b/\sqrt T$, then
\begin{equation}\label{norm:eq:thick}
 \|\Theta_K-\ee_0\|^2
 \le\left(\mathcal D_r(b)+\tau_r T^{-(r-4)/4}\right)
          T^{(r-2)/4}.
\end{equation}
\item If $S$ is even positive definite of rank $r-1$ and
$\det S\le b\sqrt T$, then
\begin{equation}\label{norm:eq:thin}
 \|\Theta_S-\ee_0\|^2
 \le\left(b^{(r-3)/2}\mathcal B_{r-1}
       +\tau_{r-1}\sqrt b\,T^{-(r-4)/4}\right)
          T^{(r-3)/4}.
\end{equation}
\end{enumerate}
In particular, for each triple $(r,b,T_0)$ below and every $T\ge T_0$,
the factors in parentheses in \eqref{norm:eq:thick} and
\eqref{norm:eq:thin} satisfy the strict upper bounds in the last
two columns:
\begin{equation}\label{norm:eq:thick-thin-table}
\begin{array}{c|c|c|c|c}
 r&b&T_0&\text{bound in \eqref{norm:eq:thick}}&
 \text{bound in \eqref{norm:eq:thin}}\\\hline
 9&8/5&200&0.017&0.014\\
 17&3/2&100&0.0015&0.0013
\end{array}.
\end{equation}
\end{corollary}
\begin{proof}
For \textup{(i)}, put $\delta=\det K$ and $u_T=b/\sqrt T$, and let
\[
 B(R)=\#\{v\in K':0<v^2\le R^2\},\qquad
 p_j=\sqrt{\beta_{r-j}/\delta},\qquad R_j=p_j^{1/j}.
\]
In Lemma~\ref{norm:lem:general}, the choice
$\alpha=(r-2)/4$ is admissible for $j\ge3$, but not for $j=2$.
We therefore refine the counting bound in dimensions at most two.

The hypothesis gives $\lambda_i(K')\ge\sqrt{u_T}$ for every $i$.
If the vectors counted by $B(R)$ span a space of dimension at most
two, \eqref{norm:eq:ball-count} gives
\[
 B(R)+1\le\left(1+\frac{2R}{\sqrt{u_T}}\right)^2,
 \qquad B(R)=0\quad(R<\sqrt{u_T}).
\]
Combining this with the higher-dimensional terms in
\eqref{r9:eq:n23}, we obtain
\[
 B(R)\le H(R)+
 \sum_{j=3}^r\frac{3^jR^j}{p_j}\mathbf1_{\{R\ge R_j\}},\qquad
 H(R)=\left(\frac{4R^2}{u_T}+\frac{4R}{\sqrt{u_T}}\right)
       \mathbf1_{\{R\ge\sqrt{u_T}\}}.
\]
Since $H$ is nondecreasing, its boundary term can be absorbed
by extending the integral to infinity. Thus its contribution is at most
\[
\begin{aligned}
 &I_r(1/4)H(1/2)
 +\int_0^{1/2}H(R)
       \left(-\frac{d}{dR}I_r(R^2)\right)\,dR\\
 \le &(r-2)C_r
 \int_{\sqrt{u_T}}^\infty
 \left(\frac{4R^2}{u_T}+\frac{4R}{\sqrt{u_T}}\right)
 R^{1-r}\,dR=
 C_r u_T^{-(r-2)/2}
 \left(\frac{4(r-2)}{r-4}+\frac{4(r-2)}{r-3}\right).
\end{aligned}
\]
Combining this with the $j\ge3$ terms and the tail bound
$\tau_r\sqrt\delta\le\tau_r\sqrt T$ gives
\[
 \|\Theta_K-\ee_0\|^2
 \le\left(\mathcal D_r(b)+\tau_rT^{-(r-4)/4}\right)
       T^{(r-2)/4}.
\]

For \textup{(ii)}, apply \eqref{norm:eq:det-only} in rank $r-1$:
\[
 \|\Theta_S-\ee_0\|^2
 \le\mathcal B_{r-1}(\det S)^{(r-3)/2}
       +\tau_{r-1}\sqrt{\det S}
 \le b^{(r-3)/2}\mathcal B_{r-1}T^{(r-3)/4}
       +\tau_{r-1}\sqrt b\,T^{1/4}.
\]
Factoring out $T^{(r-3)/4}$ gives \eqref{norm:eq:thin}.

Both factors in parentheses decrease with $T$, so it suffices to
verify \eqref{norm:eq:thick-thin-table} at $T=T_0$ in each row.
Using \eqref{r17:eq:tail}, these endpoint inequalities
are verified by rational interval evaluation.
\end{proof}

We now combine these norm bounds with the coefficient estimates
of Section~\ref{r9:sec:modular} to bound the Fourier coefficients
of $\Theta_L$. Put
\[
 \mathcal N(L)=2\zeta(r-2)C_r\lambda_1(L')^{2-r}
             +\sum_{j=2}^r C_{r,j,(r-2)/2}D^{(r-2)/2}+\tau_r\sqrt D.
\]
The constants are those of
Lemma~\ref{norm:lem:general}. One may also use the determinant-only
bound in \eqref{norm:eq:det-only}, or any sharper bound for the theta norm.
For a compatible positive index, define
\[
 \mathcal E_L(t)=
 \left(\sqrt{\mathcal N(L)}+c_\kappa\right)
 \frac{(4\pi t)^{(\kappa-1)/2}}{\sqrt{\Gamma(\kappa-1)}}
 \left(1+C(\kappa)\frac{\mathcal H(A_L)t}{\sqrt D}\right)^{1/2}.
\]

\begin{proposition}\label{gen:prop:theta-coefficients}
For every compatible positive index $(\mu,t)$,
\begin{equation}\label{gen:eq:theta-coefficients}
 \left|a_L(\mu,t)-a_{E_L}(\mu,t)\right|
 \le\mathcal E_L(t).
\end{equation}
If $L\otimes\Z_p$ contains a unimodular summand of rank $2m$ at
every prime, even at $2$, for some integer $m\ge3$, then
\begin{equation}\label{gen:eq:theta-coefficients-lower}
 a_L(\mu,t)\ge
 \max\left\{0,
 \frac{(2\pi)^\kappa\eta(m)}{\Gamma(\kappa)\sqrt D}
         t^{\kappa-1}-\mathcal E_L(t)\right\}.
\end{equation}
For $r=8,9,16,17$, the upper bounds in
\eqref{r17:eq:cuspcoef} may replace $C(\kappa)$ in the error term.
\end{proposition}
\begin{proof}
By
Lemma~\ref{r9:lem:cusp-projection} and
Lemma~\ref{norm:lem:general},
\[
 \|F_L\|\le\sqrt{\mathcal N(L)}+c_\kappa.
\]
Apply Lemma~\ref{r9:lem:poincare-bound} to $F_L$; then
Lemma~\ref{r9:lem:local-lower} together with $a_L(\mu,t)\ge0$ yields
\eqref{gen:eq:theta-coefficients-lower}.
\end{proof}

For an arbitrary $L$, a local density may vanish. For example, if
$L=\langle4\rangle^{\oplus r}$, then $Q(L)\subset2\Z$ and
$a_L(0,t)=0$ for every positive odd integer $t$, although the index
is compatible. Thus, no uniform positive lower bound is possible
without additional local hypotheses.

\subsection{Cyclic discriminant groups}\label{sec:lorentzian-estimates}

For the root-decomposition problem, cyclicity provides the 
local hypotheses required. The least corank-one determinant then selects one of the two norm bounds above. We give the argument without assuming a Lorentzian embedding.

Let $n\ge1$ be an integer and let $K$ be an even positive definite
lattice of rank $8n+1$ with $A_K$ cyclic. Put $\delta=\det K$, let $T\ge\delta$,
and choose a primitive corank-one sublattice $S\subset K$ of least
determinant $c$. Remark~\ref{r9:rem:local-hypotheses} gives
$\ell(A_S)\le2$.
Thus $K$ and $S$ have unimodular local summands of ranks $8n$ and
$8n-2$, respectively, even at $2$. The local-density and
discriminant-group estimates apply with
\[
 \begin{array}{c|cc}
 L&m&\mathcal H(A_L)<\\\hline
 K&4n&(9/8)\delta^{1/6}\\
 S&4n-1&(9/8)c^{1/4}
 \end{array}.
\]
These properties hold for every primitive corank-one $S$.
For any $b>0$, the minimizing choice and \eqref{r9:eq:n6} give
\begin{equation}\label{analytic:eq:minimality-bridge}
 \begin{cases}
 c\ge b\sqrt T\ \Longrightarrow\
       \lambda_1(K')^2=c/\delta\ge b/\sqrt T,\\
 c\le b\sqrt T\ \Longrightarrow\
       \det S\le b\sqrt T.
 \end{cases}
\end{equation}
Hence one of the two bounds in Corollary~\ref{norm:cor:thick-thin}
always applies. If $c>1$, Lemma~\ref{r9:lem:det-bounds} gives $c\ge4$,
since $\rk S=8n$. For every compatible index $(\mu,t)$ with
$t\ge T/8+1$, Lemma~\ref{r9:lem:fibre} provides a fiber index
\[
 t-Q(w)\ge t-\frac{\delta}{8c}\ge\frac{3T}{32}+1.
\]
These are the hypotheses used in the coefficient comparison below.

\begin{proposition}\label{gen:prop:analytic-range}
Let $K$ be an even positive definite lattice of rank $8n+1$, where $n\ge1$ is an integer, and suppose that $A_K$ is cyclic. Suppose that the least determinant $c$ of a
primitive corank-one sublattice is greater than one. For $T\ge\det K$
and every compatible index $(\mu,t)$ with $t\ge T/8+1$,
\begin{equation}\label{analytic:eq:general-cyclic}
 a_K(\mu,t)\gg_n T^{4n-5/4}\qquad(T\longrightarrow\infty),
\end{equation}
where both the implied constant and the lower bound for $T$
depend only on $n$, not on $K$, $\mu$ or $t$.
For $n=1,2$, the explicit bounds are
\begin{equation}\label{analytic:eq:explicit-cyclic}
 a_K(\mu,t)>
 \begin{cases}
 T^{11/4}/200>2640,&n=1,\ T\ge200,\\
 T^{27/4}/10^6>62880,&n=2,\ T\ge100.
 \end{cases}
\end{equation}
\end{proposition}
\begin{proof}
Fix $b>0$. Apply \eqref{analytic:eq:minimality-bridge}: take $L=K$
if $c\ge b\sqrt T$ and $L=S$ otherwise. Put $r=\rk L$,
$\nu=r/2-1$, and
\[
 (\mu_L,t_L)=
 \begin{cases}(\mu,t),&L=K,\\(\mu_S,t-Q(w)),&L=S.\end{cases}
\]
Then $a_K(\mu,t)\ge a_L(\mu_L,t_L)$. The determinant and index
bounds established above are
\[
 \det L\le D_r(T),\qquad t_L>t_r(T),\qquad
 (D_r(T),t_r(T),\theta_r)=
 \begin{cases}
 (T,T/8,1/6),&r=8n+1,\\
 (b\sqrt T,3T/32,1/4),&r=8n.
 \end{cases}
\]
Moreover, $\mathcal H(A_L)<(9/8)(\det L)^{\theta_r}$.
Let $\eta_r>0$ be a lower bound for the local density product
supplied by Lemma~\ref{r9:lem:local-lower}. The Eisenstein and cusp
coefficient estimates give
\begin{equation}\label{analytic:eq:relative-bound}
 \left|\frac{a_{F_L}(\mu_L,t_L)}{a_{E_L}(\mu_L,t_L)}\right|^2
 \le\frac{\nu^2\Gamma(\nu)}{4\pi^{\nu+2}\eta_r^2}
 \|F_L\|^2 t_L^{-\nu}
 \left(\det L+\frac98 C(r/2)t_L
                   (\det L)^{\theta_r+1/2}\right).
\end{equation}
With $\|F_L\|$ fixed, the right side increases with $\det L$ and
decreases with $t_L$, since $\nu>1$. Corollary~\ref{norm:cor:thick-thin} and
Lemma~\ref{r9:lem:cusp-projection} give
$\|F_L\|\ll_{n,b}T^{(r-2)/8}$ for $T\ge1$. Substituting the
determinant and index bounds in \eqref{analytic:eq:relative-bound}
yields
\begin{equation}\label{analytic:eq:relative-exponents}
 \left|\frac{a_{F_L}(\mu_L,t_L)}{a_{E_L}(\mu_L,t_L)}\right|^2
 \ll_{n,b}
 \begin{cases}
 T^{5/4-2n}+T^{23/12-2n},&L=K,\\
 T^{1-2n}+T^{15/8-2n},&L=S.
 \end{cases}
\end{equation}
All four exponents are negative for $n\ge1$. Apply this construction
with $b=1$.
For sufficiently large $T$, uniformly in the stated data,
$a_L(\mu_L,t_L)\ge\tfrac12a_{E_L}(\mu_L,t_L)$.
The Eisenstein lower bounds therefore give
\[
 a_K(\mu,t)\ge a_L(\mu_L,t_L)\gg_n
 \begin{cases}T^{4n-1},&L=K,\\T^{4n-5/4},&L=S.\end{cases}
\]
This proves \eqref{analytic:eq:general-cyclic}.

For the explicit bounds, apply the same construction with
$(b,T_0)=(8/5,200)$ when $n=1$ and $(b,T_0)=(3/2,100)$ when $n=2$.
For $T\ge T_0$, Corollary~\ref{norm:cor:thick-thin} gives
\[
 \|F_L\|\le N_r(T):=\sqrt{B_r}\,T^{(r-2)/8}+\frac1{\sqrt{24}},
\]
with the following choices of $B_r$ and $\eta_r$:
\begin{equation}\label{analytic:eq:parameter-table}
 \begin{array}{c|c|c|c}
 r&B_r&\eta_r&\sqrt{\Psi_r(T_0)}<\\\hline
 8&0.014&\eta(3)&0.942\\
 9&0.017&\eta(4)&0.946\\
 16&0.0013&49/50&0.561\\
 17&0.0015&49/50&0.530
 \end{array}.
\end{equation}
By \eqref{analytic:eq:relative-bound}, the squared relative error is
at most
\begin{equation}\label{analytic:eq:relative-majorant}
 \Psi_r(T):=
 \frac{\nu^2\Gamma(\nu)}{4\pi^{\nu+2}\eta_r^2}
 N_r(T)^2t_r(T)^{-\nu}
 \left(D_r(T)+\frac98 C(r/2)t_r(T)
                         D_r(T)^{\theta_r+1/2}\right).
\end{equation}
Since
\[
 T^{-(r-2)/4}N_r(T)^2
 =\left(\sqrt{B_r}
       +\frac{T^{-(r-2)/8}}{\sqrt{24}}\right)^2
\]
decreases with $T$, and the remaining powers of $T$ are the
negative exponents in \eqref{analytic:eq:relative-exponents},
$\Psi_r(T)$ also decreases.
Using \eqref{r17:eq:cuspcoef}, the endpoint bounds in
\eqref{analytic:eq:parameter-table} therefore give, for $T\ge T_0$,
\[
 \left|\frac{a_{F_L}(\mu_L,t_L)}{a_{E_L}(\mu_L,t_L)}\right|
 \le\sqrt{\Psi_r(T)}
 \le\sqrt{\Psi_r(T_0)}<\frac{19}{20}.
\]
Consequently, \eqref{r9:eq:eisenstein-product} gives
\begin{equation}\label{analytic:eq:two-quantitative-branches}
 a_K(\mu,t)
 \ge a_L(\mu_L,t_L)
 >\frac1{20}a_{E_L}(\mu_L,t_L)
 \ge
 \frac{(2\pi)^{r/2}\eta_r}
      {20\Gamma(r/2)\sqrt{D_r(T)}}\,t_r(T)^{r/2-1}.
\end{equation}
Substituting the parameters above and using
$T^{4n-1}\ge T^{4n-5/4}$ for $T\ge1$ yields
\eqref{analytic:eq:explicit-cyclic}.
The numerical comparisons are verified by rational interval evaluation.
\end{proof}

For $\gamma=ky\in\II_{8n+1,1}$ with $\gamma^2=-T<0$,
Lemma~\ref{r9:lem:gluing-shell} supplies all the hypotheses on $K$:
$\rk K=8n+1$, $A_K$ is cyclic, $\det K=\delta\le T$, and
$R_n(\gamma)=a_K(\mu,T/8+1)$. Hence, when $c>1$,
\eqref{analytic:eq:general-cyclic} and
\eqref{analytic:eq:explicit-cyclic} give the same bounds for
$R_n(\gamma)$. The order-two property of this particular $\mu$
is not needed in the coefficient estimates.

\begin{remark}\label{analytic:rem:nonminimal}
For an arbitrary primitive corank-one sublattice $S_0\subset K$,
Lemma~\ref{r9:lem:fibre} gives only
$\lambda_1(K')^2\le\det S_0/\delta$. Therefore, a large value of $\det S_0$ does not imply the hypothesis of
Corollary~\ref{norm:cor:thick-thin} \textup{(i)}.
For example, take $y=e+df$ in $U\oplus E_8$, with $d\ge1$, and put
$h=e-df$. Then
\[
 K=E_8\oplus\Z h,\qquad h^2=2d,\qquad
 c=1,\qquad \lambda_1(K')^2=\frac1{2d}.
\]
For a root $z\in E_8$, the vector $u=z+h/(2d)\in K'$ is primitive,
since $(u,h)=1$. Its kernel $S_0=K\cap u^\perp$ satisfies
\[
 \det S_0=(\det K)u^2=4d+1.
\]
At $d=100$, this determinant exceeds $(8/5)\sqrt T$, with $T=200$,
while $R_1(y)=480$ by Corollary~\ref{gen:cor:small-split-counts}.
Thus an arbitrary choice of $S_0$ cannot replace the minimizing
choice in Proposition~\ref{gen:prop:analytic-range}.
Minimality also makes $c=1$ equivalent to the existence of a
unimodular summand of corank one; $\det S_0>1$ cannot exclude
such a summand.
\end{remark}

\section{The case of \texorpdfstring{$\II_{9,1}$}{II(9,1)}}\label{sec:ii9}
Let $M=\II_{9,1}$ and write $\gamma=ky$, where $k\ge1$, $y\in M$ is
primitive, and $y^2=-\delta<0$. Put $T=k^2\delta$,
$K=y^\perp\cap M$, and let $c$ be the least determinant of a primitive
rank-eight sublattice of $K$.
Proposition~\ref{gen:prop:analytic-range} gives $R_1(\gamma)>2640$
when $c>1$ and $T\ge200$. It remains to treat $c=1$ and the finite
range $0<T\le200$.

\subsection{Split vectors and extremal complements}
Every even positive definite unimodular lattice of rank eight is
isometric to $E_8$. Thus Proposition~\ref{gen:prop:split-vectors}
identifies $c=1$ with the existence of a decomposition
$M=U\oplus E_8$ in which $y=a_0e+b_0f$, with $a_0,b_0>0$ coprime.

\begin{lemma}\label{r9:lem:D9-count}
Let $\epsilon_1,\ldots,\epsilon_8$ be the standard coordinate vectors of $\R^8$. In the coordinate model
\[
 E_8=\{x\in\Z^8:\textstyle\sum x_i\in2\Z\}
 \ \cup\
 \{x\in(\Z+\tfrac12)^8:\textstyle\sum x_i\in2\Z\},
\]
put $\gamma_D=2e+2f+2\epsilon_1$. Then
\[
 \gamma_D^2=-4,\qquad
 \gamma_D^\perp\cap M\cong D_9,\qquad R_1(\gamma_D)=672.
\]
\end{lemma}
\begin{proof}
Write $r=pe+qf+w\in M$, where $p,q\in\Z$ and $w=(w_1,...,w_8)\in E_8$. We conclude from $(r,\gamma_D)=2(w_1-p-q)$ that $r\in\gamma_D^\perp$ if and only if $w_1=p+q$. 
Thus, $w$ has integral coordinates, and the map $r\longmapsto(p,q,w_2,\ldots,w_8)$
identifies $\gamma_D^\perp\cap M$ with the integer vectors in
$\R^9$ having even coordinate sum, preserving the squared norm.
This is $D_9$. 

For a root decomposition $\gamma_D=r+s$, the pairing equation instead
gives $w_1=p+q-1$. The norm equation becomes
\[
 (p-1)^2+(q-1)^2+w_2^2+\cdots+w_8^2=3.
\]
The parity of this equation implies the $E_8$ coordinate-sum
condition. Its solutions have exactly three coordinates equal to
$\pm1$, so their number is $2^3\binom93=672$.
\end{proof}

\subsection{Finite reduction and root counts}

Bounded-norm enumeration in the $E_{10}$ chamber is described in
\cite[\S3.1 and Theorem~3.3]{DutourHulek2024}.
We record the bounds in our coordinates.

We use the model
\begin{equation}\label{r9:eq:odd-model}
 M=\{(a;b_1,\ldots,b_{10})\in\Z^{11}:\textstyle\sum b_i=3a\},
 \quad (a;b)\cdot(a';b')=\sum b_ib_i'-aa'.
\end{equation}
It is perpendicular to the characteristic vector
$(3;1,\ldots,1)$ of square one in the odd unimodular lattice of signature $(10,1)$. In this model, write $e_i=(3;1,\ldots,0_i,\ldots,1)$ for $1\le i\le10$. 

\begin{lemma}\label{r9:lem:chamber}
Every negative-norm vector of $M$ is isometric, up to a sign change, to a vector satisfying
\begin{equation}\label{r9:eq:chamber}
 a>0,\qquad b_1\ge\cdots\ge b_{10},\qquad a\ge b_1+b_2+b_3.
\end{equation}
Such a vector has unique nonnegative integral coordinates
$c_0,\ldots,c_9$ given by
\begin{align}
 t&=3c_0+2c_1+4c_2+\sum_{j=3}^9(9-j)c_j,
 \label{r9:eq:weight-coordinates}\\
 b_i&=t+\sum_{j=i}^9c_j\quad(1\le i\le10),\qquad
 a=c_0+b_1+b_2+b_3.\nonumber
\end{align}
An empty sum is zero.
\end{lemma}
\begin{proof}
Choose the sign to make $a>0$ and sort the $b_i$ in decreasing order.
If $a<b_1+b_2+b_3$, reflection in the root
$r=(1;1,1,1,0,\ldots,0)$ changes $a$ to
\[
 a'=2a-b_1-b_2-b_3<a.
\]
Since $\sum b_i^2<a^2$,
$b_1+b_2+b_3\le\sqrt3(\sum b_i^2)^{1/2}<2a$, so $a'>0$.
Repeated sorting and reflection terminate by descent of $a$.
Now put $c_0=a-b_1-b_2-b_3$ and $c_i=b_i-b_{i+1}$ for $1\le i\le9$.
Writing $t=b_{10}$, the defining linear relation of $M$ becomes
\[
 10t+\sum_{j=1}^9jc_j
 =3c_0+9t+3c_1+6c_2+9\sum_{j=3}^9c_j,
\]
which gives \eqref{r9:eq:weight-coordinates} and uniqueness.
\end{proof}

Let $G_i$ be the vector obtained by setting $c_i=1$ and all other
coordinates equal to zero, and put $B_{ij}=-(G_i,G_j)$. These $G_i$ form a basis of $M$. In particular, a vector is primitive
if and only if $\gcd(c_0,\ldots,c_9)=1$. By direct computation, 
\begin{equation}\label{r9:eq:weight-gram}
 B_{ij}\ge0,\qquad
 (B_{00},\ldots,B_{99})=(10,4,18,42,30,20,12,6,2,0),
\end{equation}
and $(B_{09},\ldots,B_{89})=(3,2,4,6,5,4,3,2,1)$.

\begin{lemma}\label{r9:lem:enumeration}
All chamber vectors with $0<-\gamma^2\le T_0$ are enumerated by the
following recursion. Given $c_0,\ldots,c_{j-1}$, put $Q_0=C_0=0$ and for $j>0$, 
\[
 Q_j=\sum_{i,h<j}c_ic_hB_{ih},\qquad
 C_j=\sum_{i<j}c_iB_{ij}.
\]
For $j<9$, let $c_j$ range over all nonnegative integers $u$ such that
$Q_j+2uC_j+u^2B_{jj}\leq T_0$. After choosing $c_0,\ldots,c_8$,
discard if $Q_9=0$; otherwise put $\Phi=\sum_{i<9}c_iB_{i9}$ and take
\[
 0\le c_9\le\left\lfloor\frac{T_0-Q_9}{2\Phi}\right\rfloor.
\]
\end{lemma}

\begin{proof}
Each $B_{jj}$ for $j<9$ is
positive, so each range is finite and leaves out no admissible continuation. If $Q_9=0$, all preceding coordinates vanish and only isotropic multiples of $G_9$ remain. Otherwise, $\Phi>0$ and
$-\gamma^2=Q_9+2c_9\Phi$, proving the last range. 
\end{proof}

We use the following sphere identity in both signatures. Let $L$ be positive
definite even, let $\gamma=ae+bf+z\in U\oplus L$ have square $-T<0$,
and suppose $b>0$. For $q\in\Z$ and $w\in L$, the vector
$r=pe+qf+w$ gives an ordered root decomposition $(r,\gamma-r)$ if
and only if
\begin{equation}\label{comp:eq:common-sphere}
 p=\frac{T/2-aq+(z,w)}b\in\Z,\qquad
 \|bw-qz\|^2=2b^2+Tq(b-q).
\end{equation}
Indeed, the equations are $r^2=2$ and $(r,\gamma)=-T/2$;
eliminating $p$ gives \eqref{comp:eq:common-sphere}, and the same
calculation is reversible. Nonnegativity bounds $q$ to a finite interval.
The involution $r\mapsto\gamma-r$ sends $(q,w)$ to $(b-q,z-w)$.
Thus one may count $q<b/2$ with weight two and $q=b/2$ with weight one.

For the finite range, use the chamber coordinates above and
the following integral form of \eqref{comp:eq:common-sphere}.

\begin{lemma}\label{comp:ii9-shell}
Let $\gamma=(x_0;x_1,\ldots,x_{10})$ be a chamber vector with
$\gamma^2=-T<0$, and put
\[
 a=x_9,\qquad b=x_{10},\qquad
 Z_i=4x_i-x_0-a-b\quad(1\le i\le8).
\]
Then $b>0$. Ordered root decompositions of $\gamma$ are in
bijection with integers $q$, $j\in\{0,1,2,3\}$, and
$n_1,\ldots,n_8$ satisfying
\begin{align}
 \sum_{i=1}^8(4bn_i-bj-qZ_i)^2
   &=16\bigl(2b^2+Tq(b-q)\bigr),\label{comp:eq:ii9-sphere}\\
 \sum_{i=1}^8n_i&\equiv3j\pmod4,\label{comp:eq:ii9-parity}\\
 p:=\frac{T/2-aq+\frac1{16}\sum_iZ_i(4n_i-j)}b&\in\Z.
 \label{comp:eq:ii9-integrality}
\end{align}
Only $q$ with $2b^2+Tq(b-q)\ge0$ can occur.
\end{lemma}
\begin{proof}
Put $e=e_{10}$ and $f=e_9$. Then $\gamma=ae+bf+z$, with $b>0$
by \eqref{r9:eq:weight-coordinates}. The map
$y\mapsto(y_i-y_0/4)_{i=1}^8$ identifies $U^\perp$ isometrically with
\[
 \bigcup_{j=0}^3
 \left\{n-\frac j4(1,\ldots,1):n\in\Z^8,\
                    \sum_i n_i\equiv3j\pmod4\right\},
 \qquad z_i=Z_i/4.
\]
Substitute $w_i=n_i-j/4$ in \eqref{comp:eq:common-sphere}.
The four cosets are disjoint and $p$ is uniquely determined.
\end{proof}

\begin{lemma}\label{comp:ii9}
There are $3321$ chamber vectors with $0<T\le200$. Their
representation numbers satisfy
\[
\begin{array}{lrr}
\toprule
\text{class}&\text{number of vectors}&R_1(\gamma)\\
\midrule
c_0=\cdots=c_7=0,\ \gcd(c_8,c_9)=1&180&480\\
\gamma=G_1&1&672\\
\text{other primitive vectors}&2924&\ge674\\
\gamma=2G_8&1&2640\\
\text{other nonprimitive vectors}&215&\ge2642\\
\bottomrule
\end{array}
\]
\end{lemma}

\begin{proof}[Computer-assisted proof]
The program \nolinkurl{verify-ii9-extremal-480.py}, supplied with
this paper as ancillary material, implements
Lemma~\ref{r9:lem:enumeration} with $T_0=200$ and enumerates the
finite integer spheres in Lemma~\ref{comp:ii9-shell}.

Primitivity is equivalent to $\gcd(c_0,\ldots,c_9)=1$.
Since $G_8=e_9+e_{10}$ and $G_9=e_{10}$, the first row consists
of primitive split vectors.
The program generates all $3321$ inputs and checks the row sizes
in the table. For each input, it reconstructs the recorded root pairs
and verifies lattice membership, both root norms, and distinctness
of the first roots using exact integer arithmetic.

For the three exact rows, enumeration with cutoffs $481,673,2641$
returns $480,672,2640$, respectively, so the relevant spheres have
been exhausted. In the remaining primitive and nonprimitive cases,
it yields at least $673$ and $2641$ distinct first roots. The
fixed-point-free involution $r\mapsto\gamma-r$ makes $R_1(\gamma)$
even, giving the stated lower bounds.
\end{proof}

\begin{proof}[Proof of Corollary~\ref{cor:ii9-sharp}]
If $c=1$, Proposition~\ref{gen:prop:split-vectors} and
Corollary~\ref{gen:cor:small-split-counts} give
$\gamma=k(a_0e+b_0f)$ with $a_0,b_0,k>0$, $\gcd(a_0,b_0)=1$, and
\[
 R_1(\gamma)=480\quad(k=1),\qquad
 R_1(\gamma)\ge2640\quad(k\ge2),
\]
with equality in the second case precisely when
$(k,a_0,b_0)=(2,1,1)$.

If $c>1$, Proposition~\ref{gen:prop:analytic-range} gives
$R_1(\gamma)>2640$ for $T\ge200$. In the remaining range,
Lemmas~\ref{r9:lem:chamber} and~\ref{comp:ii9} give
$R_1(\gamma)\ge672$, with equality only at $G_1$ up to isometry,
and $R_1(\gamma)>2640$ for nonprimitive $\gamma$.

The only chamber vectors with $T=4$ are $G_1$ and $G_8+G_9$.
The latter has perpendicular $E_8\oplus\langle4\rangle\not\cong D_9$,
so Lemma~\ref{r9:lem:D9-count} identifies $G_1\sim\gamma_D$
and gives the stated perpendicular characterization of this orbit.
Finally, for every $T=2m>0$, the primitive split vector $e+mf$
has $R_1(e+mf)=480$.
\end{proof}

\section{The case of \texorpdfstring{$\II_{17,1}$}{II(17,1)}}\label{sec:ii17}

Let $M=\II_{17,1}=U\oplus E_8\oplus E_8$, and write
$\gamma=ky$ with $y$ primitive, $y^2=-\delta$ and $T=k^2\delta>0$.
Put $K=y^\perp\cap M$, and let $c$ be the least determinant of a
primitive rank-sixteen sublattice of $K$.
The case $c=1$ is covered by
Proposition~\ref{gen:prop:split-vectors} and
Corollary~\ref{gen:cor:small-split-counts};
for $c>1$ and $T\ge100$,
Proposition~\ref{gen:prop:analytic-range} gives $R_2(\gamma)>62880$.

\begin{lemma}\label{r17:lem:unimodular-classification}
Every even positive-definite unimodular lattice of rank sixteen is
isometric to $E_8\oplus E_8$ or $D_{16}^{+}$.
\end{lemma}
\begin{proof}
This is Witt's classification; see
\cite[Chapter~V, \S2.3, Example~\textup{(ii)}]{Serre}.
Here
$D_{16}^{+}=D_{16}\cup((\tfrac12,\ldots,\tfrac12)+D_{16})$.
The two even index-two extensions of $D_{16}$ are isometric,
since a coordinate sign change exchanges its spinor classes.
\end{proof}

\subsection{The finite range}\label{r17:sec:finite}

Use the coordinate model
\[
 E_8=\{x\in\Z^8:\textstyle\sum x_i\in2\Z\}
 \ \cup\
 \{x\in(\Z+\tfrac12)^8:\textstyle\sum x_i\in2\Z\}.
\]
For the standard coordinate vectors $\epsilon_i$, take
\begin{equation}\label{r17:eq:simple}
 \alpha_1=\tfrac12(\epsilon_1-\epsilon_2-\cdots-\epsilon_7+\epsilon_8),
 \quad \alpha_2=\epsilon_1+\epsilon_2,\quad
 \alpha_i=\epsilon_{i-1}-\epsilon_{i-2}\quad(3\le i\le8).
\end{equation}
These roots form a basis of $E_8$. Let $(\omega_i,\alpha_j)=\delta_{ij}$ and
\begin{equation}\label{r17:eq:weights}
 G=((\alpha_i,\alpha_j))^{-1}=((\omega_i,\omega_j)),\qquad
 (G_{11},\ldots,G_{88})=(4,8,14,30,20,12,6,2).
\end{equation}
The highest root is
$\theta=\sum m_i\alpha_i=\epsilon_7+\epsilon_8$, where
\[
 (m_1,\ldots,m_8)=(2,3,4,6,5,4,3,2).
\]
Primed symbols refer to the second $E_8$ factor. Put
\[
 \mathcal A_b=\left\{\sum_{i=1}^8c_i\omega_i:
 c_i\in\Z_{\ge0},\ \sum_{i=1}^8m_ic_i\le b\right\}.
\]

\begin{lemma}\label{r17:lem:reduction}
Every negative-norm vector in $M$ is isometric, up to sign, to
\begin{equation}\label{r17:eq:reduced}
 \gamma=ae+bf-w-w',\qquad a\ge b\ge1,\qquad w,w'\in\mathcal A_b.
\end{equation}
\end{lemma}
\begin{proof}
Since $2ab=w^2+(w')^2+T>0$, arrange $a,b>0$ and choose an
isometric representative minimizing $a+b$.
Finite Weyl transformations and $e\leftrightarrow f$ then make
$w,w'$ dominant and $a\ge b$.
If $t=(w,\theta)-b>0$, reflection in $e-\theta$ sends $\gamma$ to
\[
 (a-t)e+bf-(w-t\theta)-w',\qquad
 2b(a-t)=(w-t\theta)^2+(w')^2+T>0,
\]
contradicting minimality. Thus $(w,\theta)\le b$, and likewise for
$w'$, proving the assertion.
\end{proof}

\begin{lemma}\label{r17:lem:finitebound}
For $w=\sum c_i\omega_i\in\mathcal A_b$, put $D(w)=b^2-w^2$.
Then
\begin{equation}\label{r17:eq:deficit}
 D(w)\ge \frac{b(b-2c_1)}9\ge0.
\end{equation}
Writing $a=b+l$ in \eqref{r17:eq:reduced}, we have
\begin{equation}\label{r17:eq:Tfinite}
 T=D(w)+D(w')+2lb.
\end{equation}
All reduced vectors with $0<T\le T_0$ therefore occur in the finite list
\begin{equation}\label{r17:eq:enumerate}
 \begin{gathered}
 1\le b\le9T_0,\qquad
 w,w'\in\mathcal A_b,\qquad D(w)+D(w')\le T_0,\\
 0\le l\le
 \left\lfloor\frac{T_0-D(w)-D(w')}{2b}\right\rfloor,
 \qquad T>0.
 \end{gathered}
\end{equation}
The two $E_8$ factors may be interchanged.
\end{lemma}
\begin{proof}
The scaled alcove $b^{-1}\mathcal A_b$ lies in the simplex with
vertices $0,\omega_i/m_i$, whose squared norms are
\[
 0,\quad1,\quad\tfrac89,\quad\tfrac78,\quad\tfrac56,\quad
 \tfrac45,\quad\tfrac34,\quad\tfrac23,\quad\tfrac12.
\]
For $h=b-2c_1>0$, write
\[
 \frac wb=\left(1-\frac hb\right)\frac{\omega_1}{2}+\frac hb z,
 \qquad z\in\operatorname{conv}\{0,\omega_i/m_i:2\le i\le8\}.
\]
Convexity gives $w^2/b^2\le1-h/(9b)$; the case $h=0$ is immediate.
This proves \eqref{r17:eq:deficit}, and expansion gives
\eqref{r17:eq:Tfinite}.
If either $h$ or $h'=b-2c'_1$ is positive, it is an integer and
$T\ge b/9$. Otherwise $D(w)=D(w')=0$, so $l\ge1$ and $T\ge2b$.
Thus $b\le9T_0$, and all remaining ranges follow from
\eqref{r17:eq:Tfinite}.
\end{proof}

We now enumerate root decompositions using
\eqref{comp:eq:common-sphere}.

\begin{lemma}\label{comp:ii17-shell}
Let $\gamma=ae+bf+W/2\in M$, with $W\in\Z^{16}$, $b>0$ and
$\gamma^2=-T<0$. Its ordered root decompositions correspond
bijectively to $q\in\Z$ and $X\in\Z^{16}$ satisfying
\begin{align}
 \sum_{i=1}^{16}(bX_i-qW_i)^2
   &=8b^2+4Tq(b-q),\label{comp:eq:ii17-sphere}\\
 p:=\frac{2T-4aq+W\cdot X}{4b}&\in\Z.
 \label{comp:eq:ii17-integrality}
\end{align}
In each eight-coordinate block, the entries of $X$ have the same
parity and their sum is divisible by four.
The first root is $pe+qf+X/2$.
Only the finitely many $q$ for which the right side of
\eqref{comp:eq:ii17-sphere} is nonnegative occur.
\end{lemma}
\begin{proof}
Substitute $z=W/2$ and $w=X/2$ in
\eqref{comp:eq:common-sphere}.
The block conditions are exactly $X/2\in E_8\oplus E_8$,
and $p$ is uniquely determined.
\end{proof}

We now describe the bounded chamber vectors. Put
\[
 y_D=2e+2f-\omega_1-\omega'_8,\qquad
 y_A=3e+3f-\omega_2-\omega'_2,\qquad
 y_4=3e+3f-\omega_2-\omega'_7.
\]

\begin{lemma}\label{comp:ii17}\label{r17:prop:finite-counts}
For $T_0=100$, \eqref{r17:eq:enumerate} gives $74029$ inputs after
ordering the two $E_8$ factors. Their counts are
\[
\begin{array}{c|rrrr}
 R_2(\gamma)&960&1632&4488&\ge4490\\\hline
 \text{vectors}&160&2&1&73866.
\end{array}
\]
The $160$ vectors with count $960$ are primitive and lie in
an integral hyperbolic summand as in
Proposition~\ref{gen:prop:split-vectors}; their complements are
$E_8\oplus E_8$ in $80$ cases and $D_{16}^{+}$ in $80$ cases.
The two $1632$ inputs are $y_D,y_A$, both of square $-2$;
the $4488$ input is $y_4$, of square $-4$.
For the $607$ nonprimitive inputs, the counts are
\[
\begin{array}{c|rrr}
 R_2(\gamma)&55488&62880&\ge62882\\\hline
 \text{vectors}&2&2&603.
\end{array}
\]
The two inputs with count $55488$ are $2y_D$ and $2y_A$. 
\end{lemma}

\begin{proof}[Computer-assisted proof]
In the same ancillary source package,
\nolinkurl{enumerate-counts.py}, using the lattice data in
\nolinkurl{verify-ii17-finite.py}, generates all inputs in
\eqref{r17:eq:enumerate}, with the two $E_8$ factors ordered.
It produces $74029$ inputs, of which $607$ are nonprimitive;
primitivity is tested by
$\gcd(a,b,c_1,\ldots,c_8,c'_1,\ldots,c'_8)=1$.

The program \nolinkurl{count-bounds.cpp} performs the integer-sphere
counts in \eqref{comp:eq:ii17-sphere}, retaining
\eqref{comp:eq:ii17-integrality} and the block conditions.
The involution $(q,X)\mapsto(b-q,W-X)$ allows $q\le b/2$,
with weight two for $q<b/2$ and one for $q=b/2$.
Each enumerated root pair is reconstructed, and its lattice membership
and both root norms are checked by exact integer arithmetic.
Alternatively, the program uses the exact split formula of
Proposition~\ref{gen:prop:split-vectors} after verifying an
integral hyperbolic plane containing the primitive part $y$ of
$\gamma$: it checks vectors $u,v\in M$ and positive coprime
integers $A,B$ satisfying
\[
 u^2=v^2=0,\qquad (u,v)=-1,\qquad y=Au+Bv.
\]

Counting to $5000$ for all inputs and to $62882$ for nonprimitive
inputs gives the tables: values below the cutoff are exact,
whereas the cutoff is a certified lower bound.
The program \nolinkurl{verify-count-results.py} checks the output
tables and the integral hyperbolic splittings of all $160$ inputs
with count $960$, giving $80$ inputs at each cusp.
The cusp types and exceptional inputs are identified by
Lemmas~\ref{r17:lem:second-cusp} and~\ref{r17:lem:exceptional-shells}.
\end{proof}

\subsection{The two cusps and the exceptional perpendicular lattices}

For primitive isotropic $u$, put $P_u=(u^\perp\cap M)/\Z u$.
The two possibilities are $E_8\oplus E_8$ and $D_{16}^{+}$, with
$P_e\cong E_8\oplus E_8$.
Number the following roots as $\rho_1,\ldots,\rho_{19}$:
\begin{equation}\label{r17:eq:nineteen-roots}
 \alpha_1,\ldots,\alpha_8,e-\theta,\quad
 \alpha'_1,\ldots,\alpha'_8,e-\theta',\quad f-e.
\end{equation}

\begin{lemma}\label{r17:lem:second-cusp}
The primitive isotropic vector
$e_D=2e+2f-\omega_1-\omega'_1$ satisfies $P_{e_D}\cong D_{16}^{+}$.
The roots in \eqref{r17:eq:nineteen-roots} perpendicular to $e_D$,
with the first of them omitted, form a simple system of type
$D_{16}$; their reflection closure has $480$ roots.
\end{lemma}
\begin{proof}
Put $f_D=e_D+\alpha_1$. Direct calculation gives
\[
 e_D^2=f_D^2=0,\qquad(e_D,f_D)=-1,\qquad
 (e_D,\rho_i)=-\delta_{i1}-\delta_{i,10}.
\]
Thus $U_D=\Z e_D+\Z f_D$ splits off and
$P_{e_D}\cong U_D^\perp$ is unimodular.
Deleting node $2$ as well leaves a $D_{16}$ Gram matrix;
its image in $P_{e_D}$ has index two.
The spinor extension is $D_{16}^{+}$, and its added coset has
minimum squared norm $4$. Thus $P_{e_D}$ has exactly the $480$ roots of $D_{16}$.
The selected roots map to this root system.
\end{proof}

For $y\in\{y_D,y_A,y_4\}$, let $K_y=y^\perp\cap M$.
The roots in \eqref{r17:eq:nineteen-roots} perpendicular to $y$
generate a full-rank root lattice $R_y$:
\[
\begin{array}{c|c|c|c}
 y&-y^2&R_y&[K_y:R_y]\\\hline
 y_D&2&E_7\oplus D_{10}&2\\
 y_A&2&A_{17}&3\\
 y_4&4&E_6\oplus A_{11}&3
\end{array}
\]
These data follow from the Gram matrices and
\[
 (y_D,\alpha_1)=(y_A,\alpha_2)=(y_4,\alpha_2)=-1,\qquad
 \det K_y=-y^2,\qquad [K_y:R_y]^2=\frac{\det R_y}{-y^2}.
\]

\begin{lemma}\label{r17:lem:exceptional-shells}
The root systems of $K_{y_A},K_{y_D},K_{y_4}$ are precisely those
displayed above, and
\[
 R_2(y_A)=R_2(y_D)=1632,\qquad R_2(y_4)=4488,\qquad
 R_2(2y_A)=R_2(2y_D)=55488.
\]
\end{lemma}
\begin{proof}
Write $R=R_y$ and $H=K_y/R\subset A_R$.
The overlattice correspondence
\cite[Propositions~1.4.1--1.4.2]{Nikulin} gives
$q_R|_H=0$ and $A_{K_y}=H^\perp/H$.
The group $A_{K_y}$ is cyclic by Lemma~\ref{r9:lem:gluing-shell}.
Checking isotropy, up to isometries of $R$, gives the following
subgroups and preimages $\mathcal C_y$ of the nonzero order-two class:
\[
\begin{array}{c|c|c}
 R&H&\mathcal C_y\\\hline
 A_{17}&\langle6\rangle&\{3,9,15\}\\
 E_7\oplus D_{10}&\langle(\eta_7,s_+)\rangle
       &\{(0,s_-),(\eta_7,v)\}\\
 E_6\oplus A_{11}&\langle(\eta_6,4)\rangle
       &\{(0,6),(\eta_6,10),(2\eta_6,2)\}
\end{array}
\]
Here $\eta_7,\eta_6$ generate $A_{E_7},A_{E_6}$, and $v,s_+,s_-$
are the vector and spinor classes of $D_{10}$.
Every nonzero class in $H$ has minimum squared norm $4$,
so the extensions add no roots.

For $A_{N-1}$, the minimum squared norm and number of minima
in class $j\in\Z/N\Z$ are
\begin{equation}\label{r17:eq:A-cosets}
 \left(\frac{j(N-j)}N,\binom Nj\right).
\end{equation}
The corresponding data for the nonzero classes of $E_7,E_6$ and
the vector and spinor classes of $D_{10}$ are, respectively,
\[
 (3/2,56),\qquad(4/3,27),\qquad(1,20),\qquad(5/2,512).
\]
These follow by direct coordinate counting, using
$E_7=\{x\in E_8:x_7=x_8\}$,
$E_6=\{x\in E_8:x_6=x_7=x_8\}$ and the standard model of $D_{10}$.
The same models give $756$ and $3380$ square-$4$ vectors in
$E_7$ and $D_{10}$.

Lemma~\ref{r9:lem:gluing-shell} requires squared norms
$5/2,5/2,3$ in $\mathcal C_y$ for $y_A,y_D,y_4$.
The class $9$ of $A_{17}$ has minimum $9/2$ and contributes nothing;
the remaining shell counts are
\[
 \begin{aligned}
 R_2(y_A)&=2\binom{18}{3}=1632,\\
 R_2(y_D)&=512+56\cdot20=1632,\\
 R_2(y_4)&=\binom{12}{6}+2\cdot27\binom{12}{2}=4488.
 \end{aligned}
\]
For the doubles the class is zero and the squared norm is $4$,
giving
\[
 \begin{aligned}
 R_2(2y_A)&=\binom{18}{2}\binom{16}{2}+2\binom{18}{6}=55488,\\
 R_2(2y_D)&=756+3380+126\cdot180+56\cdot512=55488.
 \end{aligned}
\]
\end{proof}

\begin{proof}[Proof of Corollary~\ref{cor:ii17-sharp}]
The split formula gives
$\gamma\in\mathcal U\ \Longleftrightarrow\ (c,k)=(1,1)$.
For $c=1$,
\[
 R_2(\gamma)=960\quad(k=1),\qquad
 R_2(\gamma)\ge62880\quad(k\ge2).
\]
For $c>1$, the analytic bound excludes every claimed minimum when
$T\ge100$, while Lemma~\ref{comp:ii17} gives the bounds and equality
candidates for $0<T\le100$.
Lemma~\ref{r17:lem:exceptional-shells} identifies their perpendiculars
and verifies their counts.

The root systems $A_{17}$ and $E_7\oplus D_{10}$ distinguish
$y_A$ from $y_D$, hence also their doubles; the finite list leaves
one orbit represented by $y_4$. This proves all equality statements.
Finally $e+(T/2)f\in\mathcal U$ has count $960$ for every positive
even $T$.
\end{proof}

\section{The general case of \texorpdfstring{$\II_{8n+1,1}$}{II(8n+1,1)}}\label{sec:higher-rank}

We construct rootless unimodular summands in every rank $8n$,
$n\ge3$, and classify the exceptions of sufficiently large negative
norm in Theorem~\ref{thm:classification} \textup{(v)}.

\begin{lemma}\label{gen:lem:rootless}
For every $n\ge3$ there is a positive definite even unimodular lattice
$P_n$ of rank $8n$ containing no vectors of square $2$.
\end{lemma}
\begin{proof}
Choose even positive definite unimodular lattices $P_3,P_4,P_5$
of minimum square $4$ in ranks $24$, $32$, and $40$, respectively:
the Leech lattice \cite{Leech67} and the constructions in
\cite{KochNebe1993,Ozeki1988}; see also
\cite[\S3, p.~855, and \S4, p.~861]{Ozeki1989}.
Write $n=3a+j$ with $a\ge0$ and $j\in\{3,4,5\}$.
Then $P_n=P_3^{\oplus a}\oplus P_j$ has rank $8n$ and minimum
square at least $4$.
\end{proof}

\begin{theorem}\label{gen:thm:higher-counterexamples}
For every $n\ge3$, there is an identification
$\II_{8n+1,1}=U\oplus P_n$ such that the following statements hold:
\begin{enumerate}
\item $R_n(ke)=0$ for every nonzero integer $k$;
\item for every integer $d\ge1$, the primitive vector $e+df$ has
square $-2d$ and satisfies $R_n(e+df)=0$.
\end{enumerate}
In particular, $\mathcal R_n+\mathcal R_n\ne \II_{8n+1,1}$ for every
$n>2$.
\end{theorem}
\begin{proof}
For $P_n$ from Lemma~\ref{gen:lem:rootless}, the indefinite
classification gives $\II_{8n+1,1}\cong U\oplus P_n$
\cite[Chapter~V, \S2.2, Theorems~5--6]{Serre}.
The cusp at $e$ has quotient $P_n$, so
Proposition~\ref{gen:prop:nonnegative} gives $R_n(ke)=0$ for $k\ne0$.
Proposition~\ref{gen:prop:split-vectors} gives
\[
 R_n(e+df)=2\#\{z\in P_n:z^2=2\}=0,
 \qquad (e+df)^2=-2d.
\]
The coefficient of $e$ is one, so $e+df$ is primitive.
\end{proof}

\begin{proof}[Proof of Theorem~\ref{thm:classification}]
Parts \textup{(ii)} and \textup{(iii)} follow from
Lemma~\ref{r9:lem:gluing-shell} and
Proposition~\ref{gen:prop:nonnegative}.
For $n=1,2$, the latter proposition and
Corollaries~\ref{cor:ii9-sharp} and~\ref{cor:ii17-sharp}
give a representation of every vector.
For $n\ge3$, Theorem~\ref{gen:thm:higher-counterexamples}
gives a primitive counterexample of every negative even square and
a primitive isotropic counterexample. This proves \textup{(i)} and
\textup{(iv)}.

It remains to prove \textup{(v)}. Write $\gamma=ky$ with $y$ primitive, put $T=-\gamma^2$, and let
$c$ be the least determinant of a primitive corank-one sublattice
of $y^\perp\cap\II_{8n+1,1}$.
If $c>1$, Lemma~\ref{r9:lem:gluing-shell} and
Proposition~\ref{gen:prop:analytic-range} give
\[
 R_n(\gamma)\gg_n T^{4n-5/4}\qquad(T\to\infty),
\]
so these vectors are represented for all sufficiently large $T$.

Suppose $c=1$. Proposition~\ref{gen:prop:split-vectors} gives
$\gamma=k(a_0e+b_0f)$ in $U\oplus P$, with $a_0,b_0>0$ coprime.
If $k=1$, its exact formula gives
\[
 R_n(\gamma)=2\#\{z\in P:z^2=2\},
\]
which proves the desired criterion in the primitive case.

It remains to treat $k\ge2$. For fixed $n$ there are only finitely
many isometry classes of $P$; see \cite[Chapter~V, \S2.3]{Serre}.
The scalar theta-series estimate
\cite[Chapter~VII, \S6.5, Theorem~8 and Corollaries~1--2]{Serre}
therefore gives, uniformly in $P$,
\[
 a_P(0,m)=[q^m]\Theta_P
 =-\frac{8n}{B_{4n}}\sigma_{4n-1}(m)+O_n(m^{2n})
 \gg_n m^{4n-1}\qquad(m\to\infty),
\]
where $B_{4n}<0$ is the Bernoulli number. Here
$\sigma_{4n-1}(m)\ge m^{4n-1}$ and $4n-1>2n$.
Taking $j=\lfloor k/2\rfloor$ in the split formula gives, for sufficiently
large $T$,
\[
 m=1+a_0b_0\lfloor k^2/4\rfloor
 \ge1+\frac{k^2a_0b_0}{8}=1+\frac T{16},
 \qquad R_n(\gamma)\ge a_P(0,m)>0.\qedhere
\]
\end{proof}

Finally, we present an application. Join two lattice vectors $x,y$ when $x-y\in\mathcal R_n$.
The resulting Cayley graph has the following diameter. 

\begin{corollary}\label{gen:cor:cayley-diameter}
The Cayley graph of the additive group $\II_{8n+1,1}$ with the roots
as generators has diameter two for $n=1,2$ and diameter three for
$n\ge3$.
\end{corollary}
\begin{proof}
It suffices to show that every vector is a sum of at most three roots.
Vectors of positive square are sums of two roots by
Proposition~\ref{gen:prop:nonnegative}. For a nonzero vector
$\gamma=ae+bf+w$ of nonpositive square in $U\oplus E_8^{\oplus n}$,
we may assume $a>0$, after changing its sign and interchanging $e,f$.
Fix a root $z\in E_8^{\oplus n}$. For $t\in\Z$,
\[
 r_t=e+(t^2-1)f+tz,\qquad r_t^2=2,\qquad
 (\gamma,r_t)=-a(t^2-1)-b+t(w,z)\longrightarrow-\infty.
\]
Hence $(\gamma-r_t)^2>0$ for sufficiently large $|t|$, and
$\gamma-r_t$ is a sum of two roots. Theorem~\ref{thm:classification}
gives diameter at most two for $n=1,2$; a nonzero isotropic vector
shows that it is not one. For $n\ge3$, the nonzero isotropic
counterexamples of Theorem~\ref{gen:thm:higher-counterexamples}
require three roots.
\end{proof}

\section*{Acknowledgements}
S.Y. is supported by the National Natural Science Foundation of China (No. 12301016) and the Fundamental Research Funds for the Central Universities (No. 20720260221).

\section*{AI Use Disclosure}
The research direction and initial mathematical framework were developed by the authors.
Generative AI was used interactively in the subsequent mathematical exploration, including
the development, testing, and refinement of parts of the construction and proof. The authors
independently verified the final arguments and take full responsibility for the paper.

\end{document}